\documentclass[a4paper,12pt]{amsart}
\usepackage{amsfonts}
\usepackage{amssymb}
\usepackage{ifthen}
\usepackage{graphicx}
\newtheorem{thm}{Theorem}
\newtheorem{cor}{Corollary}
\newtheorem{lem}{Lemma}

\newtheorem{conj}{Conjecture}

\theoremstyle{definition}
\newtheorem{defn}{Definition}
\newtheorem{case}{Case}
\newtheorem{subcase}{Subcase}

\newtheorem{examp}{Example}
\newtheorem{prob}{Problem}
\newtheorem{ques}{Question}
\newtheorem{rem}{Remark}

\newcounter {own}
\def\theown {\thesection       .\arabic{own}}

\newenvironment{pf}[1][]{%
 \vskip 3mm
 \noindent
 \ifthenelse{\equal{#1}{}}%
  {{\slshape Proof. }}%
  {{\slshape #1.} }%
 }%
{\qed\bigskip}

\newcounter{alphabet}
\newcounter{tmp}
\newenvironment{Thm}[1][]{\refstepcounter{alphabet}%
\bigskip%
\noindent%
{\bf Theorem \Alph{alphabet}}%
\ifthenelse{\equal{#1}{}}{}{ (#1)}%
{\bf .} \itshape}{\vskip 8pt}

\makeatletter
\newcommand{\Ref}[1]{\@ifundefined{r@#1}{}{\setcounter{tmp}{\ref{#1}}\Alph{tmp}}}
\makeatother

\newenvironment{Lem}[1][]{\refstepcounter{alphabet}%
\bigskip%
\noindent%
{\bf Lemma \Alph{alphabet}}%
{\bf .} \itshape}{\vskip 8pt}

\newcommand{\ID}{{\mathbb D}}
\newcommand{\IZ}{{\mathbb Z}}

\def\be{\begin{equation}}
\def\ee{\end{equation}}

\newcommand{\bee}{\begin{enumerate}}
\newcommand{\eee}{\end{enumerate}}

\newcommand{\blem}{\begin{lem}}
\newcommand{\elem}{\end{lem}}
\newcommand{\bthm}{\begin{thm}}
\newcommand{\ethm}{\end{thm}}
\newcommand{\bcor}{\begin{cor}}
\newcommand{\ecor}{\end{cor}}
\newcommand{\beg}{\begin{examp}}
\newcommand{\eeg}{\end{examp}}
\newcommand{\begs}{\begin{examples}}
\newcommand{\eegs}{\end{examples}}
\newcommand{\bdefe}{\begin{defn}}
\newcommand{\edefe}{\end{defn}}
\newcommand{\bprob}{\begin{prob}}
\newcommand{\eprob}{\end{prob}}
\newcommand{\bques}{\begin{ques}}
\newcommand{\eques}{\end{ques}}
\newcommand{\bei}{\begin{itemize}}
\newcommand{\eei}{\end{itemize}}
\newcommand{\bca}{\begin{case}}
\newcommand{\eca}{\end{case}}
\newcommand{\bsca}{\begin{subcase}}
\newcommand{\esca}{\end{subcase}}

\newcommand{\bcon}{\begin{conj}}
\newcommand{\econ}{\end{conj}}
\newcommand{\bcons}{\begin{conjs}}
\newcommand{\econs}{\end{conjs}}
\newcommand{\bprop}{\begin{propo}}
\newcommand{\eprop}{\end{propo}}
\newcommand{\br}{\begin{rem}}
\newcommand{\er}{\end{rem}}
\newcommand{\brs}{\begin{rems}}
\newcommand{\ers}{\end{rems}}
\newcommand{\bo}{\begin{obser}}
\newcommand{\eo}{\end{obser}}
\newcommand{\bos}{\begin{obsers}}
\newcommand{\eos}{\end{obsers}}
\newcommand{\bpf}{\begin{pf}}
\newcommand{\epf}{\end{pf}}
\newcommand{\ba}{\begin{array}}
\newcommand{\ea}{\end{array}}
\newcommand{\beq}{\begin{eqnarray}}
\newcommand{\beqq}{\begin{eqnarray*}}
\newcommand{\eeq}{\end{eqnarray}}
\newcommand{\eeqq}{\end{eqnarray*}}

\newcommand{\ds}{\displaystyle}

\newcounter{minutes}
\divide\time by 60
\newcounter{hours}
\multiply\time by 60 \addtocounter{minutes}{-\time}

\begin{document}
\bibliographystyle{amsplain}
\title{Univalent harmonic mappings with integer or half-integer coefficients }

\thanks{
File:~\jobname .tex,
          printed: \number\day-\number\month-\number\year,
          \thehours.\ifnum\theminutes<10{0}\fi\theminutes}

\author{S. Ponnusamy  and J. Qiao $^\dagger $}
\address{S. Ponnusamy, Department of Mathematics,
Indian Institute of Technology Madras, Chennai-600 036, India.}
\email{samy@iitm.ac.in}
\address{J. Qiao, Department of Mathematics,
Hebei University, Baoding, Hebei 071002, People's Republic of China}
\email{qiaojinjing1982@yahoo.com.cn}

\keywords{Harmonic mappings, univalent, subordination, integer coefficients,
half-integer coefficients, convexity in real direction, convexity in imaginary direction\\
$
^\dagger$ {\tt Corresponding author}
}

\begin{abstract}
Let ${\mathcal S}$ denote the set of all univalent analytic functions
$f(z)=z+\sum_{n=2}^{\infty}a_n z^n$ on the unit disk $|z|<1$.
In 1946 B. Friedman found that the set $\mathcal S$ of those functions which
have integer coefficients consists of only nine functions.
In a recent paper Hiranuma and Sugawa proved that the similar set
obtained for the functions with half-integer coefficients consists of twelve
functions in addition to the nine.
In this paper, the main aim is to discuss the class of all sense-preserving univalent harmonic
mappings $f$ on the unit disk with integer or half-integer coefficients for the analytic and co-analytic
parts of $f$. Secondly, we consider the class of univalent harmonic mappings with integer coefficients,
and consider the convexity in real direction and convexity in
imaginary direction of these mappings. Thirdly, we determine the set of univalent
harmonic mappings with half-integer coefficients which are  convex in
real direction or convex in imaginary direction.
\end{abstract}


%

\maketitle
{} \hfill {\tt
}

\pagestyle{myheadings} \markboth{S. PONNUSAMY, J. QIAO}{Harmonic
mappings}

\section{Introduction}\label{sec-1}
Assume that $f=u+iv$ is a complex-valued harmonic function
defined on the unit disk $\mathbb{D}=\{z\in \mathbb{C}:\, |z|<1\}$, i.e.
$u$ and $v$ are real harmonic in $\ID$. Then $f$ admits the
decomposition $f=h+\overline{g}$, where $h$ and $g$
are analytic in $\ID$. Often $h$ and $g$ are referred to as the analytic and co-analytic
parts of $f$, respectively.
If in addition $f$ is univalent in $\ID$, then $f$ has a non-vanishing Jacobian
in $\ID$, where the Jacobian of $f$ is given by
$$J_f(z)= |f_z(z)|^2-|f_{\bar{z}}(z)|^2=|h'(z)|^2-|g'(z)|^2.
$$
We say that $f$ is sense-preserving in $\ID$
if $J_f(z)>0$ in $\ID$. Moreover, the converse is also true, see \cite{Le}.
If $f$ is sense preserving, then the complex dilatation $\omega  :=g'/h'$
is analytic in $\ID$ and maps $\mathbb{D}$ into $\mathbb{D}.$

Denote  by $\mathcal {S}_H$ the class of all univalent sense-preserving
harmonic mappings $f=h+\overline{g}$ with the power series expansions
of $h$ and $g$ about the origin are given by
\be\label{eq1.1}
h(z)=z+\sum_{n=2}^{\infty}a_n z^n~ \mbox{ and }~
g(z)=\sum_{n=1}^{\infty}b_n z^n, \quad z\in\ID,
\ee
where we write for convenience $a_0=0$ and $a_1=1$. Also, let
$$\mathcal {S}_H^0=\{f\in \mathcal {S}_H:\,f_{\bar{z}}(0)=0\},
$$
so that $\mathcal {S}_H^0\subset\mathcal {S}_H$, and let
$\mathcal{S}=\{f=h+\overline{g}\in \mathcal{S}_H^0:\, g(z)\equiv 0\}.
$
Just like the class $\mathcal {S}$ has been a central object in the study of
univalent function theory, $\mathcal {S}_H^0$ plays a vital role in the study
of harmonic univalent mappings (see \cite{CS,Du1}). The Bieberbach conjecture had been a
driving force to develop the theory of univalent functions for a long time (\cite{Go,P,Du})
and was finally solved in the affirmative by Louis de Branges in 1985. On the other hand,
the corresponding coefficient conjecture for the class $\mathcal {S}_H^0$ has not been solved
even for the second coefficient of the analytic part $h$ of $f\in \mathcal {S}_H^0$ (\cite{CS,Du1}).
We say that a harmonic function $f=h+\overline{g}$ in $\ID$ has integer coefficients if
all the Taylor coefficients $a_n$ of $h$ and $b_n$ of $g$ are {\rm
(}rational{\rm )} integers. A similar convention applies when we say
$f=h+\overline{g}$ has half-integer coefficients. In 1946, Friedman \cite{Fr} proved the
following interesting result and for a simple proof of it, we refer to \cite{Li55} (see also
\cite{Roy55}).

\begin{Thm}\label{lem3.1}
If $f\in \mathcal {S}$ has integer coefficients, then $f$ is one of the nine functions from $\mathcal {S}_{\IZ}$,
where
\be\label{eq-sint}
\mathcal {S}_{\IZ} =\left \{ z,\; \frac{z}{1\pm z},\; \frac{z}{1\pm z^2},\; \frac{z}{(1\pm z)^2},\;
\frac{z}{1\pm z+z^2}\right \}.
\ee
\end{Thm}

In \cite{Jen87}, Jenkins presented a different proof of Theorem \Ref{lem3.1} extending the results
also to functions with coefficients in an imaginary quadratic extension of the rationals,
see also \cite{Sh,To}. It is a natural question to determine all functions $f=h+\overline{g}$ in
$\mathcal {S}_H$ such that $h$ and $g$ have integer coefficients. We obtain the following
surprising result.


\bthm\label{thm3.1}
If $f=h+\overline{g}\in \mathcal {S}_H$ have integer coefficients, then $f$ is one of the nine
functions from $\mathcal {S}_{\IZ}$.
\ethm

The key ingredient in the proof of Theorem \Ref{lem3.1}  is the
Area Theorem due to Gronwall \cite{Gr}. We refer to the work of Jenkins \cite{Jen87}
for some information that led to Theorem \Ref{lem3.1} and some related ideas.
Unfortunately,  there is no corresponding area theorem for the harmonic case as in the
lines of proof of Theorem \Ref{lem3.1}. So, it becomes necessary to obtain a suitable method
to obtain a proof of  Theorem \ref{thm3.1}. In Section \ref{sec-3}, we present
a proof of Theorem \ref{thm3.1} and the proof uses Theorem \Ref{lem3.1} and
a result of Rogosinski on subordination.

A univalent harmonic function $f$ in $\ID$ is said to be convex (resp. starlike, close-to-convex)
if $f$ is univalent and maps $\ID$ onto a convex (resp. starlike with respect to the origin,
close-to-convex) domain (see \cite{Du,Du1,Go,P}).
 Observe that each $f\in \mathcal {S}_{\IZ}$ is starlike in $\ID$.

\bdefe\label{def2.1} A domain $D\subset\mathbb{C}$ is called convex
in the direction $\alpha$ $(0\leq \alpha< \pi)$ if every line
parallel to the line through $0$ and $e^{i\alpha}$ has a connected
intersection with $D$.  A univalent harmonic function $f$ in
$\mathbb{D}$ is said to be {\it convex in the direction $\alpha$} if
$f(\mathbb{D})$ is convex in the direction $\alpha$.
\edefe

Obviously, every function that is convex in the direction $\alpha$ $(0\leq
\alpha< \pi)$ is necessarily close-to-convex. Clearly, a convex function
is convex in every direction.
The class of functions convex in one direction has been studied by many mathematicians
(see, for example, \cite{Do,HengSch70,Lecko02,RZ}) as a subclass of functions introduced by
Robertson \cite{Rober36}. We denote by $\mathcal{CV}(1)$ (resp. $\mathcal{CV}(i)$)
the class of functions convex \textit{in the direction of the real axis} (resp. in the direction of the
imaginary axis). Functions in these classes are referred to as \textit{convex in real direction}
and \textit{convex in imaginary direction}, respectively.


We continue to discuss the geometric property of the  functions in $\mathcal {S}_{\IZ}$
and reformulate the following version of Theorem \ref{thm3.1}.

\bthm\label{thm3.2} Let $f\in \mathcal {S}_H$ or $f\in \mathcal {S}$
be a function with integer coefficients. Then $f\in \mathcal{CV}(1)$
if and only if $f$ is one of the eight functions from the set
$$ \mathcal {S}_{\IZ} \backslash \left \{ \frac{z}{1-z^2}\right \}
$$
and $f\in \mathcal{CV}(i)$ 
if and only if $f$ is one of the four functions from the list
$$\left \{z,\frac{z}{1\pm z},\; \frac{z}{1-z^2}\right \}.
$$
\ethm

Proof of Theorem \ref{thm3.2} is an easy consequence of looking at the image domains of $f\in \mathcal {S}_{\IZ}$.
But for the sake of completeness, we shall present its proof.

Recently, Hiranuma and Sugawa \cite{HS} determined univalent
functions with half-integer coefficients.

\begin{Thm}\label{lem4.1}\cite[Theorem 1.2]{HS} Suppose that all the
coefficients $a_n$ of a function $f$ in $\mathcal {S}$ are
half-integers. Then $f$ is one of twenty one functions from $\mathcal {S}_{\IZ}\cup \mathcal{T}$,
where $\mathcal {S}_{\IZ}$ consists of nine functions given by \eqref{eq-sint}
and $\mathcal {T}$ consisting of twelve functions given by $ \mathcal {T}= \mathcal{T}_1\cup \mathcal{T}_2$,
where
\be\label{eq-half-int1}
\mathcal{T}_1 =\left \{
z\pm\frac{z^2}{2},\; \frac{z(2\pm z)}{2(1\pm z)},\; \frac{z(2\pm
z^2)}{2(1\pm z^2)}, \frac{z(2\pm z)}{2(1- z^2)},\;\frac{z(2\pm
z)}{2(1\pm z)^2}\right \}
\ee
and
\be\label{eq-half-int2}
 \mathcal{T}_2=\{f_{+}(z), f_{-}(z)\}
\ee
with
$$  f_{+}(z)=\frac{z(2-z+z^2)}{2(1- z+z^2)}~\mbox{ and }~f_{-}(z)=\frac{z(2+z+z^2)}{2(1+z+z^2)}.
$$
\end{Thm}

In \cite{OP}, the authors pointed that each $f\in \mathcal {S}_{\IZ}$
is not only starlike in $\ID$ but is also belonging to the class $\mathcal U$
of normalized analytic functions in $\ID$ satisfying the condition
$$\left | f'(z)\left (\frac{z}{f(z)} \right )^{2}-1\right | < 1
$$
for $|z|<1$. As observed in \cite{HS}, a similar observation is not possible for
the additional twelve functions belonging to the set $\mathcal{T}$.
By a careful analysis, the authors \cite{HS} showed that every $f\in\mathcal {T}_1$
is close-to-convex in $\ID$. The two univalent functions $f_{+}(z)$ and $f_{-}(z)$ in $\mathcal{T}_2$
are neither close-to-convex nor belong to $\mathcal U$.
We remark that an analytic function that is convex in one direction is necessarily
close-to-convex, but the converse is not true. Since the class of harmonic
functions convex in real direction and the class of harmonic functions convex in imaginary direction
have special role in geometric function theory, these classes of
univalent harmonic mappings can be characterized by its analytic part
and anti-analytic part (see Lemma \Ref{lem2.1}, and Lemma \Ref{lem2.4} with $\alpha =\pi/2$).
Thus, it is natural to investigate the class of all univalent
harmonic mappings with half-integer coefficients convex in real
direction or convex in imaginary direction.

\bthm\label{thm4.1} Let $f\in \mathcal {S}_H^0(\frac{1}{2}\IZ)$, i.e. $f\in \mathcal {S}_H^0$  with
half-integer coefficients. If $f$ is convex in real direction, then
$f$ is one of the twenty one functions from $\mathcal {S}_1\cup \mathcal {T}_3\cup \mathcal {T}_4$,
where
$$\mathcal {S}_1=  \left \{z,\frac{z}{1\pm z},\; \frac{z}{1+z^2},\; \frac{z}{(1\pm z)^2},\;
\frac{z}{1\pm z+z^2}\right \},
$$
$$\mathcal {T}_3= \left \{z\pm\frac{z^2}{2}, \; \frac{z(2\pm z)}{2(1\pm z)},\; \frac{z(2+z^2)}{2(1+z^2)},
\; \frac{z(2\pm z)}{2(1\pm z)^2} \right \},
$$
and
\be\label{eq-half-shear2}
\mathcal {T}_4= \left \{{\rm Re}\left ( \frac{z}{(1\mp z)^2} \right ) +i {\rm Im}\left ( \frac{z}{1\mp z} \right ),
{\rm Re}\left ( \frac{z}{1\mp z} \right ) +i {\rm Im}\left ( \frac{z}{(1\mp z)^2} \right ),
z\pm\overline{\frac{z^2}{2}} \right \}
\ee
\ethm

We remark that, in the proof, functions in $\mathcal {T}_4$ are represented in the sequence by
$f_3(z),f_6(z),f_9(z),f_{11}(z),f_{17}(z),f_{20}(z),
$
respectively. We emphasize that there exists only six functions in $\mathcal {S}_H^0(\frac{1}{2}\IZ)\cap \mathcal{CV}(1)$
that are not conformal.

\bthm\label{thm4.2} Let $f\in \mathcal {S}_H^0(\frac{1}{2}\IZ)$. If $f$ is
convex in imaginary direction, then $f$ is one of the eleven functions from $\mathcal {T}_5\cup \mathcal {T}_6$,
where
\be\label{eq-half-shear1b}
\mathcal {T}_5=\left \{ z,\frac{z}{1\pm z},\;\frac{z}{1-z^2},\;\frac{z(2\pm z)}{2(1\pm
z)},\; \frac{z(2-z^2)}{2(1-z^2)}, \; \frac{z(2\pm z)}{2(1-z^2)}
\right \}
\ee
and
$$\mathcal {T}_6=\left \{{\rm Re}\left ( \frac{z}{1-z} \right ) +i {\rm Im}\left ( \frac{z}{(1-z)^2} \right ),\;
{\rm Re}\left ( \frac{z}{1+z} \right ) +i {\rm Im}\left ( \frac{z}{(1+z)^2} \right ) \right \}.
$$
\ethm

It is worth pointing out from Theorem \ref{thm4.2} that there exists only two functions in
$\mathcal {S}_H^0(\frac{1}{2}\IZ)\cap \mathcal{CV}(i)$ that are not conformal.

We briefly describe the organization of the paper. In Section \ref{sec-2}, we recall
necessary lemmas that are required for the proofs of Theorems \ref{thm3.1} and \ref{thm3.2}
and the proofs of these theorems will be given in Section \ref{sec-3} while the proofs
of Theorems \ref{thm4.1} and \ref{thm4.2}  will be presented in Section \ref{sec-4}.

We end the section with a conjecture.

\begin{conj}
Let $f\in \mathcal {S}_H^0(\frac{1}{2}\IZ)$. Then $f$ is one of twenty seven functions from
$\mathcal {S}_{\IZ}\cup \mathcal{T}_1\cup\mathcal {T}_2 \cup \mathcal {T}_4$
where $\mathcal {S}_{\IZ},\mathcal{T}_1, \mathcal{T}_2, \mathcal {T}_4$
are given by \eqref{eq-sint}, \eqref{eq-half-int1},
\eqref{eq-half-int2}, and \eqref{eq-half-shear2} respectively.
\end{conj}

\section*{5. Illustrations through figures}\label{sec5}

Using Mathematica, we present the images of the unit disk $\ID$ under some of these mappings.

\begin{figure}
\begin{minipage}{0.45\linewidth}
\centering
\includegraphics[height=5.5cm, width=5.5cm, scale=1]{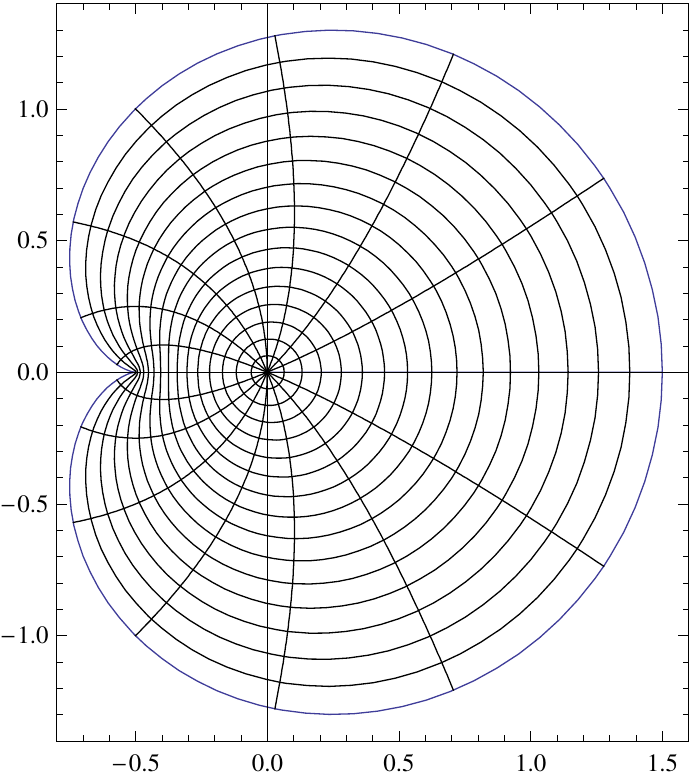}
\vspace{5pt}
$z+z^2/2$
\end{minipage}
\begin{minipage}{0.45\linewidth}
\centering
\includegraphics[height=5.5cm, width=5.5cm, scale=1]{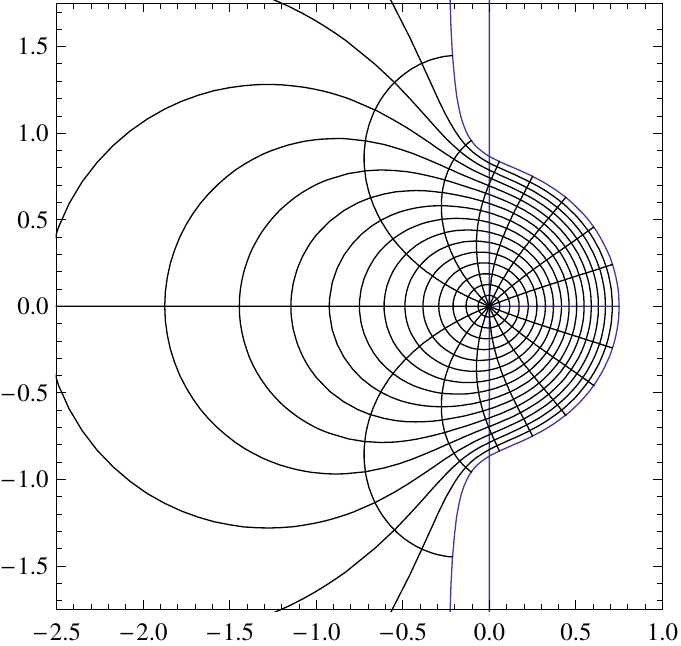}
\vspace{5pt}
$ z(2+z)/(2(1+z))$
\end{minipage}
\end{figure}

\begin{figure}
\begin{minipage}{0.45\linewidth}
\centering
\includegraphics[height=5.5cm, width=5.5cm, scale=1]{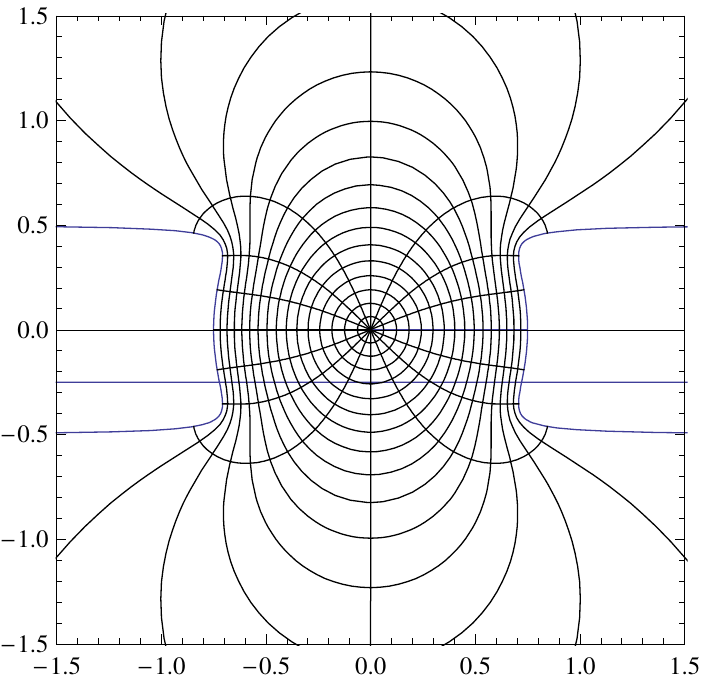}
\vspace{5pt}
$z(2+z^2)/(2(1+z^2))$
\end{minipage}
\begin{minipage}{0.45\linewidth}
\centering
\includegraphics[height=5.5cm, width=5.5cm, scale=1]{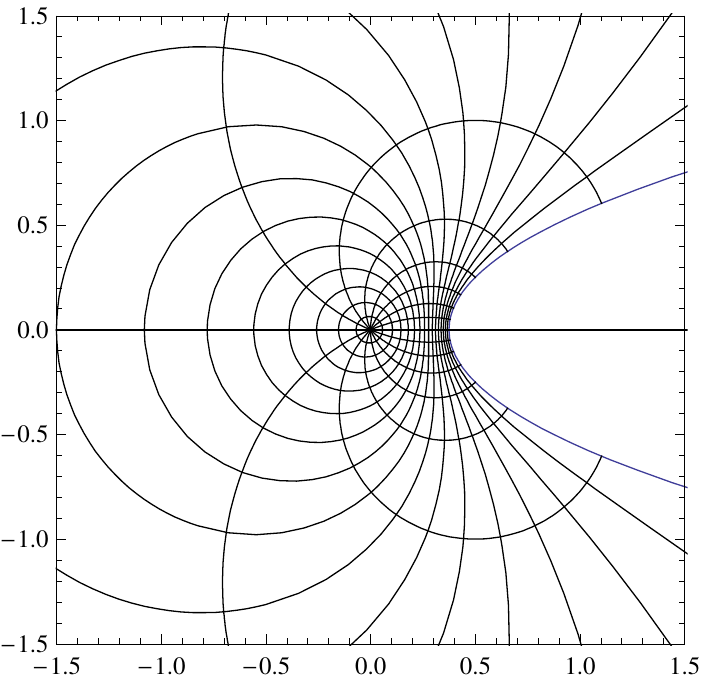}
\vspace{5pt}
$z(2+z)/(2(1+z)^2)$
\end{minipage}
\end{figure}

\begin{figure}
\begin{minipage}{0.45\linewidth}
\centering
\includegraphics[height=5.5cm, width=5.5cm, scale=1]{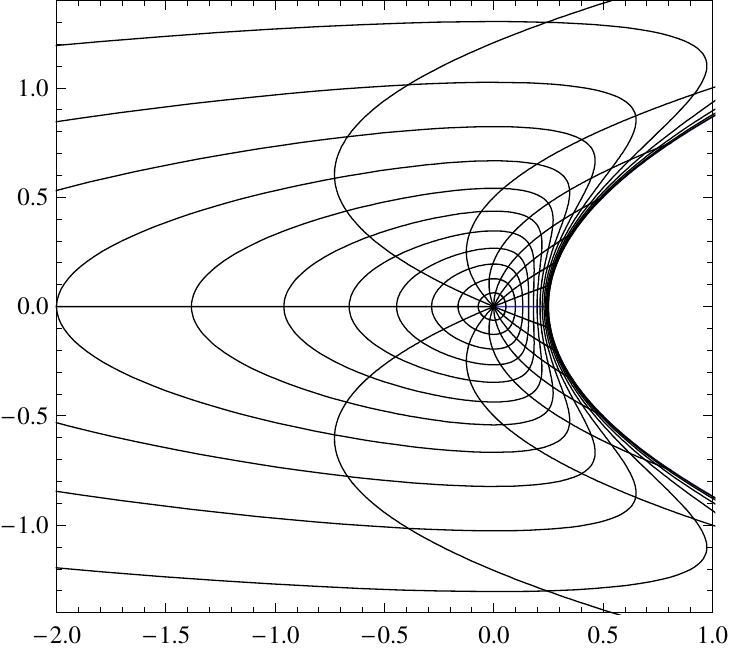}
\vspace{5pt}
${\rm Re}\,(z/(1+z)^2) + i{\rm Im}\,(z/(1+z))$
\end{minipage}
\begin{minipage}{0.45\linewidth}
\centering
\includegraphics[height=5.5cm, width=5.5cm, scale=1]{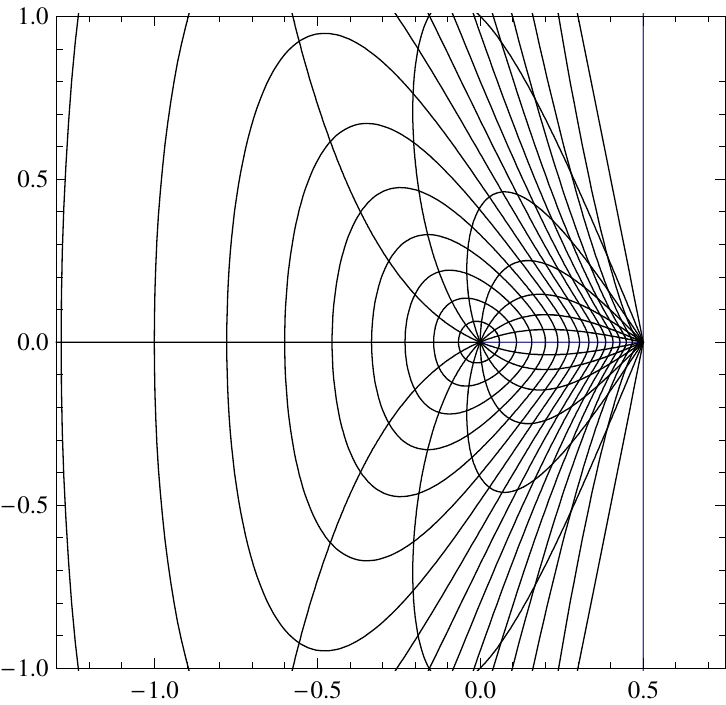}
\vspace{5pt}
${\rm Re}\,(z/(1+z)) + i{\rm Im}\,(z/(1+z)^2)$
\end{minipage}
\end{figure}

\newpage

\begin{figure}
\begin{minipage}{0.45\linewidth}
\centering
\includegraphics[height=5.5cm, width=5.5cm, scale=1]{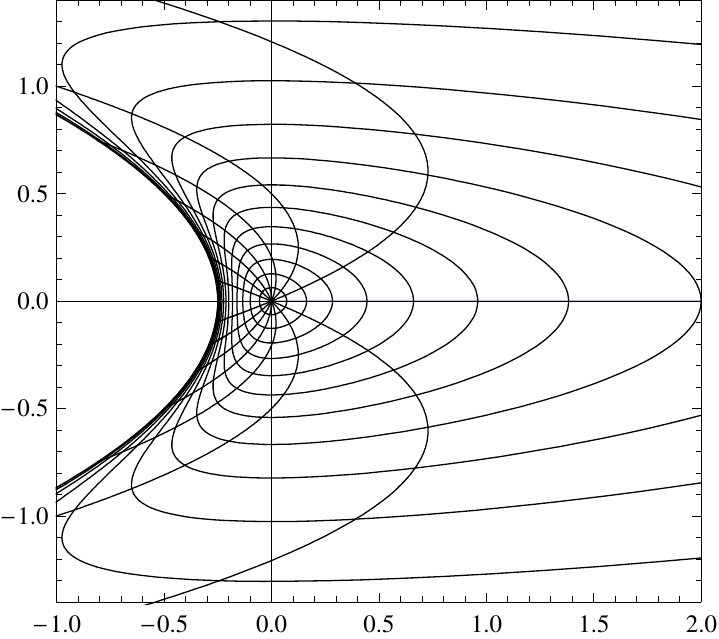}
\vspace{5pt}
${\rm Re}\,(z/(1-z)^2) + i{\rm Im}\,(z/(1-z))$
\end{minipage}
\begin{minipage}{0.45\linewidth}
\centering
\includegraphics[height=5.5cm, width=5.5cm, scale=1]{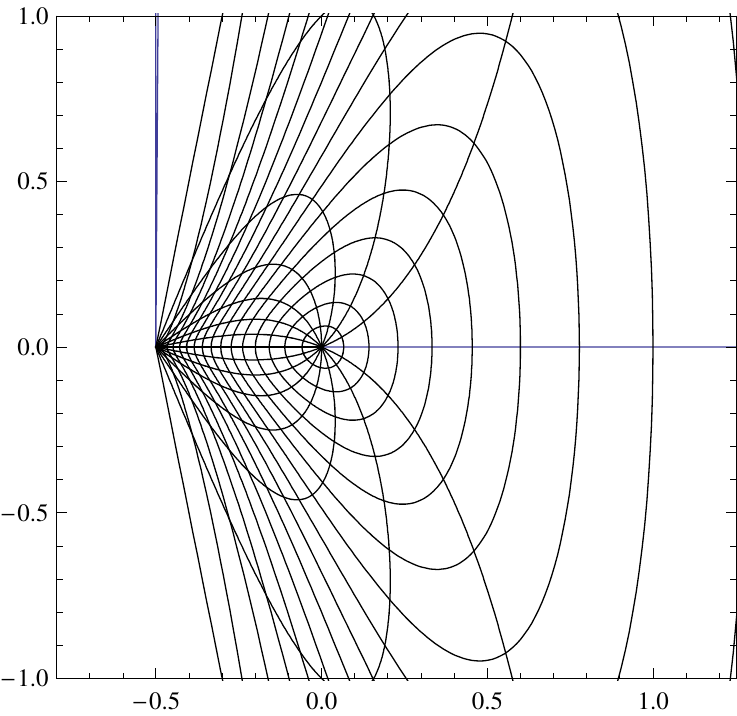}
\vspace{5pt}
${\rm Re}\,(z/(1-z)) + i{\rm Im}\,(z/(1-z)^2)$
\end{minipage}
\end{figure}

\begin{figure}
\begin{minipage}{0.45\linewidth}
\centering
\includegraphics[height=5.5cm, width=5.5cm, scale=1]{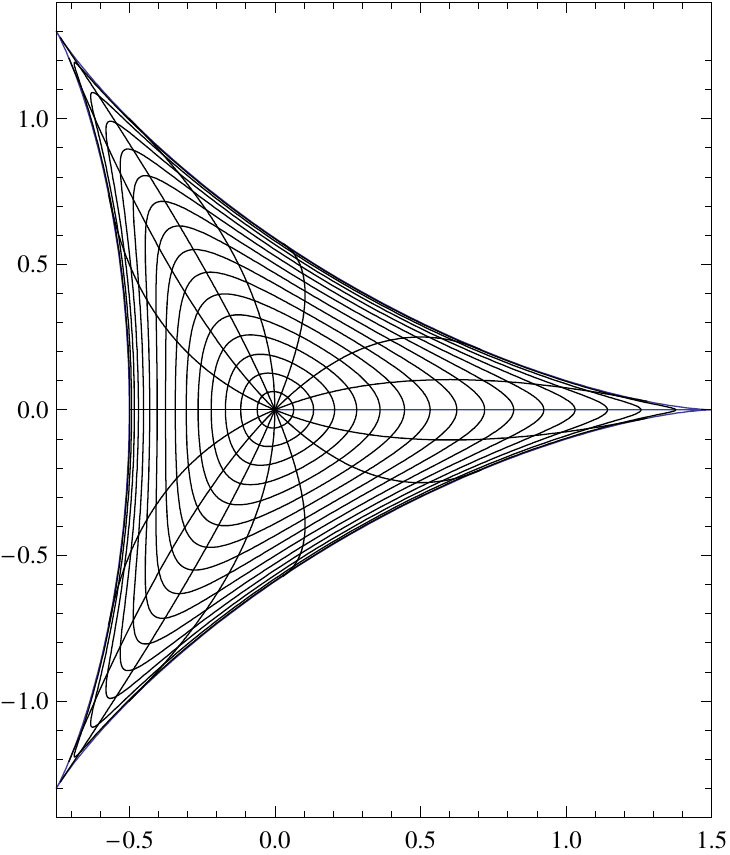}
\vspace{5pt}
$z+\overline{z^2/2}$
\end{minipage}
\begin{minipage}{0.45\linewidth}
\centering
\includegraphics[height=5.5cm, width=5.5cm, scale=1]{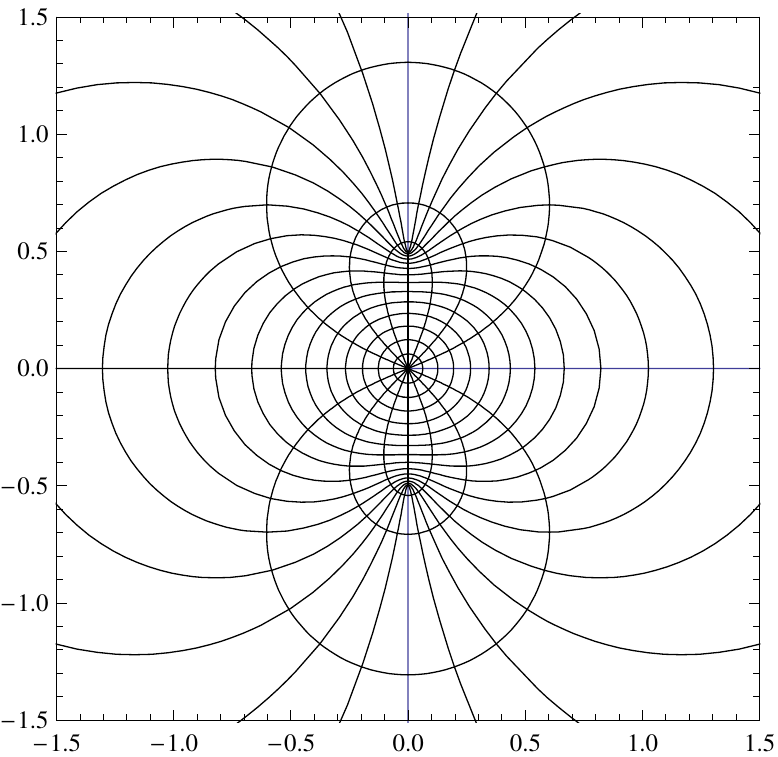}
\vspace{5pt}
$z/(1-z^2)$
\end{minipage}
\end{figure}

\begin{figure}
\begin{minipage}{0.45\linewidth}
\centering
\includegraphics[height=5.5cm, width=5.5cm, scale=1]{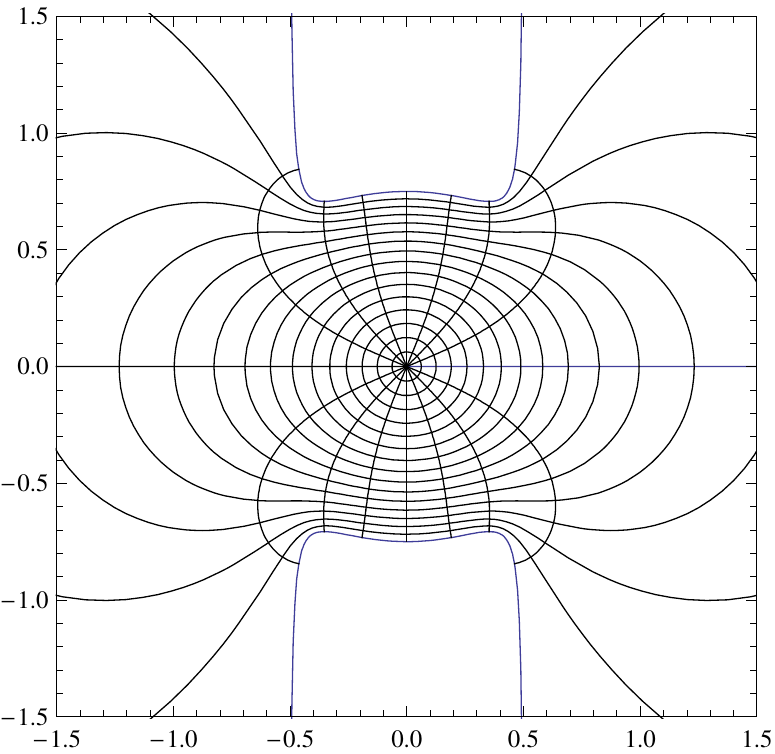}
\vspace{5pt}
$z(2-z^2)/(2(1-z^2))$
\end{minipage}
\begin{minipage}{0.45\linewidth}
\centering
\includegraphics[height=5.5cm, width=5.5cm, scale=1]{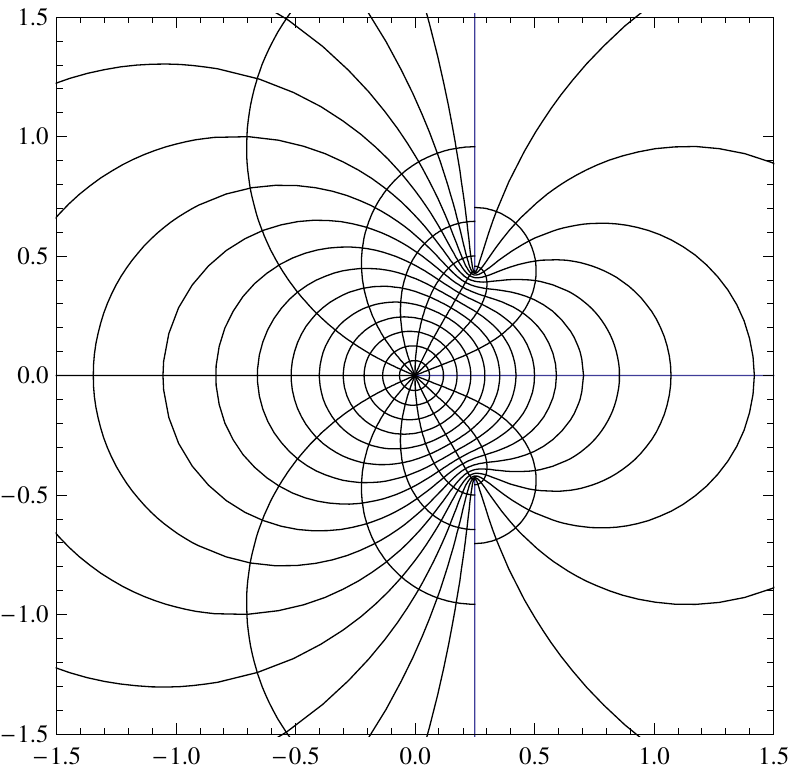}
\vspace{5pt}
$z(2-z)/(2(1-z^2))$
\end{minipage}
\end{figure}

\section{Lemmas}\label{sec-2}

%
%
%

%

%

In proving our theorems we will need a few known lemmas. The first lemma is
popularly known as Clunie and Sheil-Small's shear construction
theorem \cite[Theorem 5.3]{CS} which produces a univalent harmonic function
that maps $\ID$ onto a domain convex in the direction of the real axis.

\begin{Lem}\label{lem2.1}
A harmonic $f=h+\overline{g}$ locally univalent in $\mathbb{D}$ is a
univalent mapping of $\mathbb{D}$ onto a domain convex in the
direction of the real axis if and only if $h-g$ is a conformal
univalent mapping of $\mathbb{D}$ onto a domain convex  in the
direction of the real axis.
\end{Lem}

If $f=h+\overline{g}$ is convex in the direction $\alpha$, then
$$e^{-i\alpha}f =e^{-i\alpha}h + \overline{e^{\alpha}g} ~\mbox{ and }~
e^{-i\alpha}h-e^{i\alpha}g=e^{-i\alpha}(h-e^{2i\alpha}g)
$$
are convex in the direction of the real axis, and hence the function $h-e^{2i\alpha}g$ is convex
in the direction $\alpha$.  Thus, a natural corollary to
Lemma \Ref{lem2.1} may be stated in the following precise form so that
$f=h+\overline{g}$ is convex in the direction of the imaginary axis if and only
if $h+g$ is convex in the direction of the imaginary axis, in particular.

\begin{Lem}\label{lem2.4} 
A harmonic $f=h+\overline{g}$ locally univalent in $\mathbb{D}$ is a
univalent mapping of $\mathbb{D}$ onto a domain convex in the direction $\alpha$ $(0\leq \alpha <\pi)$
if and only if $h-e^{2i\alpha}g$ is a conformal
univalent mapping of $\mathbb{D}$ onto a domain convex in the direction $\alpha$.
\end{Lem}


Paul Greiner \cite{Gre} has constructed numerous examples using the method of shearing.
In the proofs of Theorems \ref{thm4.1} and \ref{thm4.2}, the discussion on different cases and
subcases give a number of univalent harmonic mappings that are especially convex in real
direction and/or convex in vertical direction. However, only few of them have half-integer
coefficients.

Next, we recall a useful result by Royster and Ziegler \cite{RZ} concerning analytic mappings
convex in one direction. Several particular cases of this result is helpful in
testing whether an analytic mapping is convex in a direction.

\begin{Lem}\label{lem2.2}\cite[Theroem 1]{RZ}
Let $\varphi(z)$ be a non-constant function analytic in $\mathbb{D}$. The
function $\varphi(z)$ maps  univalently $\mathbb{D}$ onto a domain convex in
the direction of imaginary axis if and only if there are numbers
$\mu$ and $\nu$, $0\leq \mu <2\pi$ and $0\leq \nu\leq \pi$, such that
\be\label{RZ-eq1}
{\rm
Re\,}\{-ie^{i\mu}(1-2ze^{-i\mu}\cos \nu+z^2e^{-2i\mu})\varphi'(z)\}\geq
0,\;\, z\in \mathbb{D}.
\ee
\end{Lem}

Again, since a function $\varphi$ is convex in real direction if and only if the function $i\varphi$ is convex
in imaginary direction, the following version is a consequence of Lemma \Ref{lem2.2}.

\begin{Lem}\label{lem2.3} Let $\varphi(z)$ be a non-constant function
regular in $\mathbb{D}$. The function $\varphi(z)$ maps  univalently
$\mathbb{D}$ onto a domain convex in the direction of real axis if
and only if there are numbers $\mu$ and $\nu$, $0\leq \mu <2\pi$ and
$0\leq \nu\leq \pi$, such that
\be\label{RZ-eq2}
{\rm
Re\,}\{e^{i\mu}(1-2ze^{-i\mu}\cos \nu+z^2e^{-2i\mu})\varphi'(z)\}\geq
0,\;\, z\in \mathbb{D}.
\ee
\end{Lem}

Using Royster and Ziegler's result,  Schaubroeck \cite{Sch2001} investigated
certain properties of the class of functions convex in a direction.

\begin{Lem}\label{lem2.5}\cite[p.87, Theorem]{Du1} For all functions
$f\in \mathcal {S}_H^0$, the sharp inequality
$|b_2|\leq \frac{1}{2}$ holds, with equality if and only if
$\omega(z)=e^{i\alpha}z$ for some real $\alpha$.
\end{Lem}

The notion of subordination is an important property in analytic
function theory, see \cite{Du,Rog}. For \textit{analytic functions} $f$
and $g$ in $\ID$, we say that $f$ is subordinate to
$g$, written $f(z)\prec g(z)$ or simply $f\prec g$, if there exists
a Schwarz' function $\varphi$ (i.e.  $\varphi$ is analytic in $\ID$
with $\varphi (0)=0$ and $|\varphi (z)|<1$ for $z\in \ID$) such that
$f(z)=g(\varphi (z)).$ The condition implies that
$f(0)=g(0)$ and $|f'(0)|\leq |g'(0)|$.
If, in addition, $g$ is univalent, then  $f\prec g$ if and
only if $f(\mathbb{D})\subset g(\mathbb{D})$ and $f(0)=g(0)$.

The following result due to Rogosinski \cite{Rog} (see also Duren \cite[p.195, Theorem 6.4]{Du})
is crucial in the proof of Theorem \ref{thm3.1}.

\begin{Lem}\label{lem2.6} If $g(z)=\sum_{n=1}^{\infty}b_nz^n$ is analytic in
$\mathbb{D}$ and $g\prec f$  for some convex function from $f\in \mathcal{S}$,
then $|b_n|\leq 1$ for $n\geq 1$.
\end{Lem}

\section{The proofs of Theorems \ref{thm3.1} and \ref{thm3.2}}\label{sec-3}

\subsection{Proof of Theorem \ref{thm3.1}} 
Let $f=h+\overline{g}\in \mathcal {S}_H$, where $h$ and $g$ have the standard form given
by \eqref{eq1.1}:
$$h(z)=z+\sum_{n=2}^{\infty}a_n z^n~ \mbox{ and }~
g(z)=\sum_{n=1}^{\infty}b_n z^n
$$
where $a_n,b_n$ are integers. Since $f$ is sense-preserving, we have  $J_f = |h'|^2-|g'|^2> 0 $
so that $\varphi'(z)\neq 0$ in $\mathbb{D}$ for $\varphi(z)=h(z)-g(z)$, and $|b_1|=|g'(0)|<|h'(0)|=1$
which implies that $b_1=0$. Thus, $f\in \mathcal {S}_H^0$. But then, by Lemma \Ref{lem2.5}, $|b_2|\leq 1/2$.
Since $b_2$ is an integer, we must have $b_2=0$.

Now, we claim that $g(z) \equiv 0$ in $\ID$. Suppose on the contrary that $g$ is not
identically zero. Because $f$ is sense-preserving, we have $|h'|=|g'+\varphi'|>|g'|$, that is
$$\left |\frac{g'}{\varphi'}+1\right |>\left |\frac{g'}{\varphi'}\right |
, ~\mbox{ i.e. }~ {\rm Re}\left \{\frac{g'(z)}{\varphi'(z)}\right \}>-\frac{1}{2} ~\mbox{ for $z\in \ID$.}
$$
In terms of subordination, we may rewrite the last inequality as
$$\frac{g'(z)}{\varphi'(z)}\prec \frac{z}{1-z} ~\mbox{ for $z\in \ID$.}
$$
Let $n_0=\min\{n:\, b_n\neq 0\}$ and observe that $\varphi'(z)=1+\sum_{n=2}^{\infty}n(a_n-b_n) z^{n-1}\neq 0$ in $\ID$
so that $1/\varphi'$ has the form
$$\frac{1}{\varphi'(z)}=1+\sum_{n=1}^{\infty}c_n z^n.
$$
Then $b_{n_0}\neq 0$ for some $n_0> 2$ and therefore, we have the representation
$$\frac{g'(z)}{\varphi'(z)}=n_0 b_{n_0}z^{n_0-1}+\sum_{n=n_0}^{\infty}d_n z^n ~\mbox{ for $z\in \ID$.}
$$
By Lemma \Ref{lem2.6}, we deduce that $|n_0b_{n_0}|\leq 1$. Since $b_{n_0}$ is an integer
and $n_0>2$, it follows that $b_{n_0}=0$ which is a contradiction. Thus, we conclude that
$g(z)\equiv 0$. Finally, the proof of the theorem follows from Theorem \Ref{lem3.1}.
\hfill $\Box$

\vspace{8pt}



\subsection{Proof of Theorem \ref{thm3.2}} 
Let $f\in \mathcal {S}_H$ (resp. $\mathcal {S}$)
with integer coefficients. Then it follows from Theorem \ref{thm3.1}
(resp. Theorem \Ref{lem3.1}) that
$$
f\in \mathcal {S}_{\IZ} = \left \{z,\; \frac{z}{1\pm z},\; \frac{z}{1\pm z^2},\; \frac{z}{(1\pm
z)^2},\; \frac{z}{1\pm z+z^2}\right \}.
$$
Clearly, $z$, $z/(1-z)$ and $z/(1+z)$ are convex in
$\mathbb{D}$. In particular, these three functions belong to $\mathcal{CV}(1)\cap \mathcal{CV}(i)$.
Next, we need to consider the remaining functions with the help of either the condition \eqref{RZ-eq1}
or \eqref{RZ-eq2}.

\begin{case} The function $\varphi(z)=z/(1+z^2)$. \end{case}
If we choose $\mu=0$ and $\nu= \pi/2$ in Lemma \Ref{lem2.3} then for this $\varphi(z)$,
we obtain that
$$ {\rm Re\,}\{(1+z^2) \varphi'(z)\}={\rm Re\,}\left \{\frac{1-z^2}{1+z^2}\right \}>0, ~z\in \ID,
$$
which implies that $\varphi(z)$ is convex in real direction. Moreover, since
$$\varphi(e^{i\theta})= \frac{e^{i\theta}}{1+e^{2i\theta}}=\frac{1}{2 \cos\theta},
$$
it follows that $\varphi(z)$ maps $\mathbb{D}$ onto the complex
plane slit along the two half-lines $y=0,~ |x|\geq 1/2$. Note that
$\varphi(z)$ is not convex in imaginary direction.

Next, we consider the case  $\psi (z)=z/(1-z^2)$ and observe that $\psi (z)=-i\varphi(iz)$.
As consequence of this observation, we see that  $\psi$ is convex in imaginary
direction and $\psi(z)$ maps $\mathbb{D}$ onto the complex plane slit along the two
half-lines $x=0, ~|y|\geq 1/2$. In particular, $\psi (z)$ is not convex in real direction.

\bca The Koebe function $k(z)=z/(1-z)^2$.\eca
It is well-known that $k(z)$ maps $\mathbb{D}$ onto the
complex plane slit along the half-line $y=0,~ x<-1/4$ and so, $k\in \mathcal{CV}(1)\backslash \mathcal{CV}(i)$.

Next, we let $p(z)=z/(1+z)^2$ and observe that $p(z)=-k(-z)$. Thus,
$p(z)$ maps $\mathbb{D}$ onto the complex plane slit
along the half-line $y=0,~ x>1/4$ and so, has similar geometric properties as that of
$k(z)$.

\bca The function $\varphi(z)=z/(1-z+z^2)$. \eca  If we choose $\mu=0$ and
$\nu=2\pi/3$ in Lemma \Ref{lem2.3}, then we have
$${\rm Re\,}\,\{(1-z+z^2)\varphi'(z) \}
={\rm Re\,}\left \{\frac{1-z^2}{1-z+z^2}\right \}>0, ~z\in \ID,
$$
because the inversion
$$\frac{1-z+z^2}{1-z^2} =\frac{1}{1+z} +\frac{z^2}{1-z^2}
$$
obviously has positive real part. Thus, by Lemma \Ref{lem2.3}, $\varphi(z)$ is convex in real direction.
Moreover, as
$$\varphi(e^{i\theta})=\frac{e^{i\theta}}{1-e^{i\theta}+e^{2i\theta}}=\frac{1}{2 \cos\theta-1},
$$
it follows easily that $\varphi(z)$ maps $\mathbb{D}$ onto the
complex plane slit along the two half-lines $y=0, ~x<-1/3$
and $y=0,\, x>1$. In particular, 
$\varphi(z)\in \mathcal{CV}(1)\backslash \mathcal{CV}(i)$.

Finally, if $\psi (z)= z/(1+z+z^2)$, then $\psi (z)= -\varphi(-z)$
so that $\psi (z)$ maps $\mathbb{D}$ onto the complex plane slit
along the two half-lines $y=0,~ x>1/3$, and $y=0,~ x<-1$. Consequently,
$\psi(z)\in \mathcal{CV}(1)\backslash \mathcal{CV}(i)$.
\hfill $\Box$


\section{The proofs of Theorems \ref{thm4.1} and \ref{thm4.2}}\label{sec-4}

We need the following lemma which gives the complete information on the set of all
univalent functions with half-integer coefficients, that are
either convex in real direction or convex in imaginary direction.

\blem\label{lem4.2}
Let $f\in \mathcal {S}$ with half-integer
coefficients. Then $f$ is convex in real direction if and only if
$f\in \mathcal {S}_1\cup \mathcal {T}_3$, where
$$\mathcal {S}_1= \mathcal {S}_{\IZ} \backslash \left \{ \frac{z}{1-z^2}\right \}
~\mbox{ and }~ \mathcal {T}_3= \mathcal {T}_{1} \backslash \left \{\frac{z(2-z^2)}{2(1-z^2)} \right \}.
$$
Moreover, $f$ is convex in imaginary direction if and only if $f$ is one of the nine functions from
$\mathcal {T}_5$, where $\mathcal {S}_{\IZ}$, $\mathcal{T}_1$ and
$\mathcal {T}_5$ are given by \eqref{eq-sint}, \eqref{eq-half-int1}, and \eqref{eq-half-shear1b},
respectively.
\elem\bpf
Assume that $f\in \mathcal {S}$ with half-integer coefficients. According to Theorem \Ref{lem4.1},
$f$ must belong to $\mathcal {S}_{\IZ}\cup \mathcal{T}_1\cup \mathcal{T}_2$,
where $\mathcal{T}_2=\{f_{+}(z), f_{-}(z)\}$ is given by  \eqref{eq-half-int2} and observe that
$$f_{+}(z)=\frac{z(2-z+z^2)}{2(1- z+z^2)}~\mbox{ and }~f_{-}(z)=\frac{z(2+z+z^2)}{2(1+z+z^2)}
$$
are not close-to-convex and hence, can neither be convex in real direction nor be convex in imaginary
direction. Thus, in view of Theorem \ref{thm3.2}, is suffices to deal with only the functions
in $\mathcal{T}_1$. Consequently, we need to check which functions in $\mathcal{T}_1$
are convex in real direction or convex in imaginary direction.

\setcounter{case}{0}

\bca The function $\varphi _1(z)=z- z^2/2.$  \eca

Choose $\mu=0$ and $\nu= 2\pi/3$ in Lemma \Ref{lem2.3} with $\varphi =\varphi _1$.
Then $\varphi _1(z)$ maps $\ID$ univalently onto a domain convex in real direction, since
$$
{\rm Re\,}\{(1+z+z^2)\varphi _1'(z)\} = {\rm Re\,}\{1-z^3\}>0.
$$
The line $x=1/2$ intersects the boundary of the image
domain $\varphi _1(\mathbb{D})$ at least with three points: $(\frac{1}{2}, 0)$, $(\frac{1}{2}, 1)$,
$(\frac{1}{2}, -1)$. It follows that $\varphi _1(z)$ is not convex in imaginary direction.

Setting $\psi _1(z)=-\varphi _1(-z)$ shows that $\psi _1(z)$ is convex in real direction
but not convex in imaginary direction.

\bca The function $\varphi _2(z)= z(2-z)/ (2(1-z))$. \eca

If we set $\mu=0$ and $\nu=0$ in Lemma \Ref{lem2.3} with $\varphi = \varphi _2$, then the condition
\eqref{RZ-eq2} is satisfied, since
$$
{\rm Re\,}\{(1-z)^2\varphi _2'(z) \}={\rm Re}\left \{(1-z)^2\frac{(1-z)^2+1}{2(1-z)^2}\right \}
=\frac{1}{2}{\rm Re\,}\{1+(1-z)^2\}>0,
$$
showing that $\varphi _2(z)$ is convex in real direction.

Similarly, if we choose $\mu=\frac{\pi}{2}$ and $\nu=\frac{\pi}{2}$ in Lemma \Ref{lem2.2} with $\varphi=
\varphi _2 $, then \eqref{RZ-eq1} is also satisfied, since
$${\rm Re\,}\{(1-z^2)\varphi _2'(z)\}={\rm
Re}\left \{(1-z^2)\frac{(1-z)^2+1}{2(1-z)^2}\right \}= \frac{1}{2}{\rm
Re}\left \{1-z^2+\frac{1+z}{1-z}\right \}>0,
$$
which gives that $\varphi _2(z)$ is also convex in imaginary direction.

Now, for the function  $\psi _2(z)= z(2+z)/ (2(1+z))$, we observe that
$\psi _2(z)=-\varphi _2(-z)$ and therefore,  $\psi _2(z)$ is convex both in real
and imaginary direction.

\bca The function $\varphi _3(z)=z(2-z^2)/(2(1-z^2))$. \eca

Choose $\mu=\frac{\pi}{2}$ and $\nu=\frac{\pi}{2}$ in Lemma \Ref{lem2.2}
with $\varphi = \varphi _3$. Then, the condition \eqref{RZ-eq1} is also satisfied, since
$$
{\rm Re\,}\{(1-z^2)\varphi _3'(z) \}={\rm Re}\left \{(1-z^2)\frac{2-z^2+z^4}{2(1-z^2)^2}\right \}
=\frac{1}{2}\,{\rm Re}\left \{1-z^2+\frac{1+z^2}{1-z^2}\right \}>0.
$$
Thus,  $\varphi _3(z)$ is convex in imaginary direction. Next, we observe that
$$ \varphi _3(e^{i\theta})=\frac{e^{i\theta}(2-e^{2i\theta})}{2(1-e^{2i\theta})}
=\frac{\cos \theta}{2}+i\left (\frac{\sin \theta}{2}+\frac{1}{4\sin \theta}\right ),
$$
and so, it is a simple exercise to verify that the line $y=3/4$  meets the
boundary of the image domain $\varphi _3(\mathbb{D})$ at least with three points:
$(\frac{\sqrt{3}}{4}, \frac{3}{4})$, $(-\frac{\sqrt{3}}{4},
\frac{3}{4})$, $(0, \frac{3}{4})$. Consequently,  $\varphi _3(z)$
is not convex in real direction.

Again if we let $\psi _3(z)= z(2+z^2)/(2(1+z^2))$, then, since  $\psi _3(z)=-i\varphi _3(iz)$,
it follows easily that $\psi _3(z)$ is convex in real
direction but is not convex in imaginary direction.

\bca The function $\varphi _4(z)=z(2-z)/(2(1-z^2))$.\eca
For this function, we see that
$$4\varphi _4(z)-1=\frac{-1+4z-z^2}{1-z^2},
$$
and
$$\varphi _4(e^{i\theta})=\frac{e^{i\theta}(2-e^{i\theta})}{2(1-e^{2i\theta})}
=\frac{2-e^{i\theta}}{-4i\sin\theta }
= \frac{1}{4}-i\frac{2-\cos \theta}{4\sin \theta}.
$$
Therefore, it follows easily that  $\varphi _4(z)$ maps the unit disk $\mathbb{D}$
onto the complex plane slit along the two half-lines $\frac{1}{4}+iy$, $|y|\geq
\sqrt{3}/4$. In particular, $\varphi _4(z)$ is convex in imaginary direction, but not convex in
real direction.

Now, for the function $\psi _4(z)=z(2+z)/(2(1-z^2))$, we see that $\psi _4(z)=-\varphi _4(-z)$
and therefore, we conclude that $\psi _4(z)$ maps $\mathbb{D}$ onto the complex plane
slit along the two half-lines $-\frac{1}{4}+iy$, $|y|\geq \sqrt{3}/4$.
In particular, $\psi _4(z)$ is convex in imaginary direction, not convex in real direction.

\bca  The function $\varphi _5(z)=z(2-z)/(2(1-z)^2)=(1-(1-z)^2)/(2(1-z)^2)$. \eca
If we choose $\mu=0$ and $\nu=0$ in Lemma \Ref{lem2.3} with $\varphi=\varphi _5$,
then
$${\rm Re\,}\{(1-z)^2\varphi _5'(z)\}= {\rm Re }\left \{\frac{1}{1-z}\right \}>0,
$$
which shows that $\varphi _5(z)$ is convex in real direction. Further, it
is easy to check that the boundary of the image domain $\varphi _5(\mathbb{D})$
is the parabola $8u+16v^2+3=0$. In particular, $\varphi _5(z)$ is not convex in imaginary
direction.

Finally, if we let $\psi _5(z)=z(2+z)/(2(1+z)^2)$, then we have $\psi _5(z)=-\varphi _5(-z)$
and therefore, $\psi _5(z)$ is convex in real direction, but is clearly not convex in imaginary direction.

The desired conclusion follows if we compile all the above cases together.
\epf

Now, we are ready to prove Theorems \ref{thm4.1} and \ref{thm4.2}.

\subsection{Proof of Theorem \ref{thm4.1}}
Let $f=h+\overline{g}\in \mathcal {S}_H^0(\frac{1}{2}\IZ)$ such that $f\in \mathcal{CV}(1)$,
and further let $\varphi=h-g$, where $h$ and $g$ have the standard form given
by \eqref{eq1.1} with $b_1=0$:
$$h(z) =z+\sum_{n=2}^{\infty}a_n z^n ~\mbox{ and }~g(z)=\sum_{n=2}^{\infty}b_n z^n.
$$
According to Lemma \Ref{lem2.1}, since $f$ is univalent,
$f\in \mathcal{CV}(1)$ if and only if $\varphi$ is a conformal univalent mapping such
that $\varphi \in \mathcal{CV}(1)$.
Moreover, as $f\in \mathcal {S}_H^0(\frac{1}{2}\IZ)$, we observe  that the analytic
function $\varphi=h-g$ also has half-integer coefficients. It follows from Lemma
\Ref{lem4.2} that $\varphi \in \mathcal {S}_1\cup \mathcal {T}_3$, where $\mathcal {S}_1$
and $\mathcal {T}_3$ are defined as in Lemma \ref{lem4.2}.
As in the proof of Theorem \ref{thm3.1}, since $f$ is sense-preserving, an analysis shows that
$$\frac{g'(z)}{\varphi'(z)}\prec \frac{z}{1-z}, \quad z\in \ID,
$$
where $\varphi(z)$ is locally univalent analytic function in $\ID$ satisfying the normalization
condition $\varphi(0)=0$ and $\varphi'(0)=1$. By Lemma \Ref{lem2.6}, we obtain that
$|b_2|\leq 1/2$. As $2b_2$ is an integer, we must have either $b_2=0$ or $b_2=\pm 1/2$.
These two conditions are necessary for $f$ to belong to $\mathcal {S}_H^0(\frac{1}{2}\IZ)$.

\setcounter{case}{0}

\bca  The case $b_2=0$.\eca

We claim that $g(z) \equiv 0$. Suppose on contrary
that $g(z_0)\neq 0$ for some $z_0\in\mathbb{D}$. Then, we may let
$n_0=\min\{n: b_n\neq 0\}$ so that
$$\frac{1}{\varphi'(z)}=1+\sum_{n=1}^{\infty}c_n z^n.
$$
Note that $b_{n_0}\neq 0$ $(n_0\geq 3)$ and
$$\frac{g'(z)}{\varphi'(z)}=n_0
b_{n_0}z^{n_0-1}+\sum_{n=n_0}^{\infty}d_n z^n.
$$
By Lemma \Ref{lem2.6}, we have $|n_0b_{n_0}|<1$, since $b_{n_0}$ is a
half-integer and $n_0\geq 3$, it follows that $b_{n_0}=0$ which is a
contradiction. Thus, $g(z) \equiv 0$ in $\ID$ and hence,  by Lemma \Ref{lem4.2},
the analytic part $h$ of $f$ must belong to one of the functions from the
set $\mathcal {S}_1\cup \mathcal {T}_3$.

The conditions $b_2=\pm 1/2$ are necessary for $f$ to belong to $\mathcal {S}_H^0(\frac{1}{2}\IZ)$,
but not sufficient to claim $f\in \mathcal {S}_H^0(\frac{1}{2}\IZ)$ as we can see below in a number of
examples.

\bca The case $b_2=\pm 1/2$.\eca

Since $b_2=\pm 1/2$, by Lemma \Ref{lem2.5},  we deduce that $\omega(z)=e^{i\alpha}z$. Because
$2a_n$ and $2b_n$ are integers and the dilatation $\omega(z)$ satisfies the condition
$g'(z)=\omega(z) h'(z)$, $\alpha$ is either $0$ or $\pi$ and thus, $\omega(z)=\pm z$.

Solving $g'(z)=\omega(z) h'(z)$ together with $\varphi(z)=h(z)-g(z)$ gives us the harmonic function $f(z)$
in a convenient form:
\be\label{eq-half-shear1}
f(z)=h(z)+\overline{g(z)}=2{\rm Re\,}h(z)-\overline{\varphi (z)},
~\mbox{ with }~h(z)=\int_0^z\frac{\varphi '(t)}{1-\omega (t)}\, dt
\ee
where $\omega(z)=\pm z$. We divide the remaining part of the proof into several subcases and in all these
cases, we need first to compute the integral in \eqref{eq-half-shear1}.

\setcounter{case}{0}

\bsca The case $\varphi(z)=z$ and $\omega(z)=\pm z$.  \esca

Evaluating the integral in \eqref{eq-half-shear1} with $\varphi(z)=z$ and $\omega(z)=z$ yields
$$h(z)=-\log (1-z)
$$
and when $\varphi(z)=z$ and $\omega(z)=-z$, it gives $h(z)=\log (1+z)$. Thus,
the corresponding functions in questions in these two cases are given by
$$f_1(z)=-2\log |1-z| -\overline{z}
~\mbox{ and }~
f_2(z)=2\log |1+z| -\overline{z}.
$$
Note that both the analytic and co-analytic parts of the univalent harmonic functions $f_1(z)$ and $f_2(z)$
do not have half-integer coefficients, although they are convex in real direction. It can be easily seen that
$f_1(z)$ maps $\ID$ onto a domain bounded by three concave arcs, with
cusps at the points $-\log 2 \pm i$ and $\infty$. 


\bsca The case $\varphi(z)=z/(1-z)$ and $\omega(z)=\pm z$. \esca

The analytic part $h(z)$ in \eqref{eq-half-shear1} when $\varphi(z)=z/(1-z)$ and $\omega(z)=\pm z$
takes the form
\beqq
h(z)&=&\int_0^z\frac{dt}{(1-t)^2(1-\omega (t))} \\
&=& \left\{ \begin{array}{rl}
\ds \frac{1}{2}\left ( \frac{1}{(1-z)^2} -1 \right ) & \mbox{for $\omega(z)=z$},  \\
\ds \frac{z}{2(1-z)} +\frac{1}{4}\log \left (\frac{1+z}{1-z}\right ) & \mbox{for $\omega(z)=-z$}.
\end{array}\right.
\eeqq
In the case when $\omega(z)=z$,  a simplification with the above $h$ gives
$$f_3(z)=h(z)+\overline{g(z)})= \frac{1}{2} \left ( \frac{z}{(1-z)^2} +\frac{z}{1-z} \right ) +
\frac{1}{2}\overline{ \left ( \frac{z}{(1-z)^2}- \frac{z}{1-z}\right ) }
$$
so that the function
$$f_3(z)= \sum_{n=1}^{\infty}\left (\frac{n+1}{2}\right ) z^n +
\overline{\sum_{n=2}^{\infty}\left (\frac{n-1}{2}\right ) z^n}
$$
has half-integer coefficients. Thus, $f_3\in \mathcal {S}_H^0(\frac{1}{2}\IZ)$ and
the function $f_3(z)$ may be equivalently written as
$$f_3(z)= {\rm Re}\left ( \frac{z}{(1-z)^2} \right ) +i {\rm Im}\left ( \frac{z}{1-z} \right ).
$$
If we let $\ell (z)=z/(1-z)$ and $k(z)=z/(1-z)^2$, we see that
$$f_3(z)=\frac{1}{2} \left ( \ell (z) + k(z)\right ) -
\frac{1}{2}\overline{\left ( \ell (z)- k(z)\right )}=
{\rm Re}\,k(z) +i {\rm Im}\,\ell (z).
$$
which is indeed not convex (see the relevant figure in Section \ref{sec5}). It is worth remarking that although
the function
$$L(z)= {\rm Re}\,\ell (z) +i {\rm Im}\,k(z)=
\sum_{n=1}^{\infty}\left (\frac{n+1}{2}\right ) z^n -\overline{\sum_{n=2}^{\infty}\left (\frac{n-1}{2}\right ) z^n}
$$
is known to be a well-known extremal function for the coefficient inequality
for the class of convex functions from $\mathcal {S}_H^0$ (see \cite{CS,Du1}). However
coefficients of $f_3(z)$ do satisfy the necessary coefficient conditions for convex functions in
$\mathcal {S}_H^0$ without $f_3(z)$ being convex in $\ID$.

In the case when $\omega(z)=-z$, it is clear that the Taylor coefficients of the
corresponding $h$ above does not have half-integer coefficients and so, the corresponding harmonic
function
$$f_4(z)= \frac{1}{2}\log \left |\frac{1+z}{1-z}\right | +i {\rm Im}\left ( \frac{z}{1-z} \right )
$$
does not belong to  $\mathcal {S}_H^0(\frac{1}{2}\IZ)$.


\bsca The case $\varphi(z)=z/(1+z)$ and $\omega(z)=\pm z$. \esca
The analytic part $h(z)$ in \eqref{eq-half-shear1} when $\varphi(z)=z/(1+z)$ takes the form
\beqq
h(z)&=&\int_0^z\frac{dt}{(1+t)^2(1-\omega (t))} \quad  (\omega(z)=\pm z)\\
&=& \left\{ \begin{array}{rl}
\ds \frac{z}{2(1+z)} +\frac{1}{4}\log \left (\frac{1+z}{1-z}\right ) & \mbox{for $\omega(z)=z$}\\
\ds \frac{1}{2}\left (\frac{z}{1+z}+ \frac{z}{(1+z)^2} \right ) & \mbox{for $\omega(z)=-z$} .
\end{array}\right.
\eeqq
By a computation, we can easily see that the corresponding harmonic functions are given by
$$f_5(z)= \frac{1}{2}\log \left |\frac{1+z}{1-z}\right | +i {\rm Im}\left ( \frac{z}{1+z} \right )
$$
and
$$f_6(z)=h(z)+\overline{g(z)}= \frac{1}{2} \left ( \frac{z}{(1+z)^2} +\frac{z}{1+z} \right ) +
\frac{1}{2}\overline{ \left ( \frac{z}{(1+z)^2}- \frac{z}{1+z}\right )}
$$
or equivalently,
$$f_6(z)= {\rm Re}\left ( \frac{z}{(1+z)^2} \right ) +i {\rm Im}\left ( \frac{z}{1+z} \right ),
$$
respectively.
We see that $f_6(z)\in \mathcal {S}_H^0(\frac{1}{2}\IZ)$, but $f_5z)\not\in \mathcal {S}_H^0(\frac{1}{2}\IZ)$.
Observe that $f_6(z)=-f_3(-z)$ and so, $f_6(z)$ is not convex although it does satisfy
the necessary coefficient conditions for convex functions in $\mathcal {S}_H^0$


\bsca The case $\varphi(z)=z/(1+z^2)$ and $\omega(z)=\pm z$.  \esca\
The function $h(z)$ in \eqref{eq-half-shear1} in this case takes the form
\be\label{eq-half-shear1a}
h(z)=\int_0^z\frac{1-t^2}{(1+t^2)^2(1-\omega (t))}\,dt \quad  (\omega(z)=\pm z).
\ee
For $\omega(z)=z$, we write the integrand as
$$\frac{1+t}{(1+t^2)^2}=\frac{1}{4} \left [\frac{1}{(1+it)^2}+\frac{1}{(1-it)^2}+
\frac{1}{1+it}+\frac{1}{1-it} \right ] +\frac{t}{(1+t^2)^2}
$$
and the integration leads
$$h(z)=\frac{z^2}{2(1+z^2)} +\frac{z}{2(1+z^2)}+\frac{1}{4i}\log \left (\frac{1+iz}{1-iz}\right ).
$$
This gives
$$f_7(z)={\rm Re}\left ( \frac{z^2}{(1+z)^2} \right ) +\frac{1}{2}\arg \left (\frac{1+iz}{1-iz}\right ) +i
{\rm Im}\left ( \frac{z}{1+z^2} \right )
$$
and it clear that $f_7(z)\not\in \mathcal {S}_H^0(\frac{1}{2}\IZ)$, because $h$ does not have
half-integer coefficients.

In the case of $\omega(z)=-z$, simplifying the integrand in \eqref{eq-half-shear1a}, we obtain that
$$h(z)=-\frac{z^2}{2(1+z^2)} +\frac{z}{2(1+z^2)}+\frac{1}{4i}\log \left (\frac{1+iz}{1-iz}\right )
$$
and therefore, the corresponding harmonic function is
$$f_8(z)=-{\rm Re}\left ( \frac{z^2}{(1+z)^2} \right ) +\frac{1}{2}\arg \left (\frac{1+iz}{1-iz}\right ) +i
{\rm Im}\left ( \frac{z}{1+z^2} \right )
$$
which is again not in $\mathcal {S}_H^0(\frac{1}{2}\IZ)$.

\bsca The case $\varphi(z)=z/(1-z)^2$ and $\omega(z)=\pm z$. \esca
If $\omega(z)=z$, then as above we end up with
$$f_9(z)=\frac{z-\frac{1}{2}z^2+\frac{1}{6}z^3}{(1-z)^3}+\overline{\frac{\frac{1}{2}z^2+\frac{1}{6}z^3}{(1-z)^3}}.
$$
which is indeed the well-known harmonic Koebe function (with dilatation $\omega(z)=z$) and the
function has no half-integer coefficients.

If $\omega(z)=-z$, then a calculation gives
$$ h(z)=\frac{z(2-z)}{2(1-z)^2}=\frac{1}{2} \left ( \frac{z}{(1-z)^2}+\frac{z}{1-z}\right )
$$
and
$$ g(z)=h(z)-\varphi(z)=-\frac{1}{2} \frac{z^2}{(1-z)^2} =-\frac{1}{2} \left ( \frac{z}{(1-z)^2}- \frac{z}{1-z}\right )
$$
which leads to
$$f_9(z)= {\rm Re}\left ( \frac{z}{1-z} \right ) +i {\rm Im}\left ( \frac{z}{(1-z)^2} \right ),
$$
or equivalently in power series,
$$f_9(z)= \sum_{n=1}^{\infty}\left (\frac{n+1}{2}\right ) z^n -
\overline{\sum_{n=2}^{\infty}\left (\frac{n-1}{2}\right ) z^n}.
$$
Clearly, the function $f_9(z)$ belongs to $\mathcal {S}_H^0(\frac{1}{2}\IZ)$.

\bsca The case $\varphi(z)=z/(1+z)^2$ and $\omega(z)=\pm z$. \esca
As in the previous case, it follows easily that for $\omega(z)=z$ we get
$$ f_{11}(z)=\frac{z(2+z)}{2(1+z)^2}+\overline{\frac{z^2}{2(1+z)^2}} =
{\rm Re}\left ( \frac{z}{1+z} \right ) +i {\rm Im}\left ( \frac{z}{(1+z)^2} \right )
$$
and for $\omega(z)=-z$ we obtain
$$ f_{12}(z)=\frac{z+\frac{1}{2}z^2+\frac{1}{6}z^3}{(1+z)^3}+\overline{\frac{\frac{1}{2}z^2-\frac{1}{6}z^3}{(1+z)^3}}.
$$
Again, by a minor calculation, we see that $f_{11}(z)$ belongs to $\mathcal {S}_H^0(\frac{1}{2}\IZ)$
whereas $f_{12}(z)$ does not belong to $\mathcal {S}_H^0(\frac{1}{2}\IZ)$.


\bsca The case $\varphi(z)=z/(1-z+z^2)$ and $\omega(z)=\pm z$.\esca
In this case, the analytic part $h(z)$ of the corresponding harmonic mappings satisfies
$$h'(z)=\left \{ \begin{array}{rl}
\ds \frac{1+z}{(1-z+z^2)^2} & \mbox{for } \omega(z)=z,\\
\ds \frac{1-z}{(1-z+z^2)^2} & \mbox{for }  \omega(z)=-z.
\end{array}\right.
$$
It is easy to check that, for $h(z)=\sum_{n=1}^{\infty}a_kz^k$,
$$a_4=\left \{ \begin{array}{rl}
\ds -\frac{1}{4} & \mbox{for } \omega(z)=z,\\
\ds -\frac{3}{4} & \mbox{for }  \omega(z)=-z,
\end{array}\right.
$$
which implies that in these two cases the corresponding harmonic mappings
$f_{13}$ and $f_{14}$ do not have half-integer coefficients.

\bsca  The case $\varphi(z)=z/(1+z+z^2)$ and $\omega(z)=\pm z$.\esca
We see that the analytic part $h(z)$ satisfies
$$h'(z)=\left \{ \begin{array}{rl}
\ds \frac{1+z}{(1+z+z^2)^2} & \mbox{for } \omega(z)=z,\\
\ds \frac{1-z}{(1+z+z^2)^2} & \mbox{for } \omega(z)=-z.
\end{array}\right.
$$
As in the previous case,  we obtain that, for $h(z)=\sum_{n=1}^{\infty}a_kz^k$,
$$a_4=\left \{\begin{array}{rl}
\ds \frac{3}{4} & \mbox{for } \omega(z)=z,\\
\ds \frac{1}{4} & \mbox{for } \omega(z)=-z,
\end{array}\right.
$$
which shows that the corresponding harmonic mappings
$f_{15}$ and $f_{16}$ do not have half-integer coefficients.


\bsca The case $\varphi(z)=z-z^2/2$ and $\omega(z)=\pm z$.
\esca
For $\omega(z)=z$, it is easy to obtain that
$$f_{17}(z)=z+\overline{\frac{z^2}{2}}
$$
which clearly belongs to $\mathcal {S}_H^0(\frac{1}{2}\IZ)$. On the other hand,
in the case $\omega(z)=-z$, the analytic part of the corresponding
harmonic mapping $f_{18}(z)$ turns out to be
$$h(z)=2\log(1+z)-z=z-\sum_{n=2}^{\infty}\frac{(-1)^n}{k}z^k,
$$
and so  $f_{18}(z)=2{\rm Re\,} h(z) -\overline{(z-z^2/2)}$ does not have half-integer coefficients.

\bsca The case $\varphi(z)=z+z^2/2$ and $\omega(z)=\pm z$.
\esca
If $\omega(z)=z$, a computation gives $f_{19}(z)$ with the corresponding
analytic part as
$$h(z)=-2\log(1-z)-z=z+\sum_{n=1}^{\infty}\frac{1}{n}z^n.
$$
It follows that  $f_{19}(z)=2{\rm Re\,} h(z) -\overline{(z+z^2/2)}$
does not have half-integer coefficients. If $\omega(z)=-z$, then it is easy to obtain
$$f_{20}(z)=z-\overline{\frac{z^2}{2}}
$$
which is clearly in $\mathcal {S}_H^0(\frac{1}{2}\IZ)$.


\bsca The case $\varphi(z)=z(2-z)/(2(1-z))$ and $\omega(z)=\pm z$.
\esca
In these two cases, the corresponding analytic parts are obtained from
$$h'(z)=\left\{ \begin{array}{rl}
\ds \frac{(1-z)^2+1}{2(1-z)^3}  & \mbox{for } \omega(z)=z,\\
\ds \frac{(1-z)^2+1}{2(1-z)^2(1+z)}& \mbox{for } \omega(z)=-z.
\end{array}\right.
$$
Indeed integrating the last for each case  gives
$$h(z)=\left\{ \begin{array}{rl}
\ds  -\frac{1}{2}\log(1-z)+\frac{1}{4(1-z)^2}-\frac{1}{4}
  & \mbox{for } \omega(z)=z,\\
\ds \frac{5}{8}\log(1+z)-\frac{1}{8}\log(1-z)+\frac{1}{4(1-z)}-\frac{1}{4}
 & \mbox{for } \omega(z)=-z;
\end{array}\right.
$$
or equivalently,
$$h(z)=\left\{ \begin{array}{rl}
\ds  \sum_{n=1}^{\infty}\left(\frac{1}{2n}+\frac{n+1}{4}\right )z^n & \mbox{for } \omega(z)=z,\\
\ds\sum_{n=1}^{\infty}\left (\frac{(-1)^{n+1}5}{8n}+\frac{1}{8n}+\frac{1}{4}\right )z^n & \mbox{for } \omega(z)=-z.
\end{array}\right.
$$
We see that in these two cases the corresponding harmonic mappings $f_{21}(z)$ and $f_{22}(z)$
do not have half-integer coefficients.


\bsca The case $\varphi(z)=z(2+z)/(2(1+z))$ and $\omega(z)=\pm z$.
\esca
As in the previous case,  we obtain

$$
h'(z)=\left\{ \begin{array}{rl}
\ds \frac{(1+z)^2+1}{2(1+z)^2(1-z)} & \mbox{for } \omega(z)=z,\\
\ds \frac{(1+z)^2+1}{2(1+z)^3}& \mbox{for }  \omega(z)=-z.
\end{array}\right.
$$
Then if $\omega(z)=z$, we obtain
$$
h(z)=\frac{1}{8}\log(1+z)-\frac{5}{8}\log(1-z) +\frac{z}{4(1+z)}
=\sum_{n=1}^{\infty}\left (\frac{(-1)^{n+1}}{8n}+\frac{5}{8n}-\frac{(-1)^n}{4}\right )z^n
$$
and for $\omega(z)=-z$, we have
$$h(z)=\frac{1}{2}\log(1+z)-\frac{1}{4}\left (\frac{1}{(1+z)^2}-1\right )
=\sum_{n=1}^{\infty}(-1)^{n+1}\left (\frac{1}{2n}+\frac{n+1}{4}\right )z^n.
$$
Thus, in both cases the corresponding harmonic mappings $f_{23}(z)$ and $f_{24}(z)$
do not belong to $\mathcal {S}_H^0(\frac{1}{2}\mathbb{Z})$.


\bsca The case $\varphi(z)=z(2+z^2)/(2(1+z^2))$ and $\omega(z)=\pm z$. \esca
In this case, we obtain that

$$h'(z)=\left\{ \begin{array}{rl}
\ds \frac{2+z^2+z^4}{2(1+z^2)^2(1-z)}& \mbox{if } \omega(z)=z,\\
\ds \frac{2+z^2+z^4}{2(1+z^2)^2(1+z)}& \mbox{if } \omega(z)=-z.
\end{array}\right.
$$
If  $h(z)=\sum_{n=1}^{\infty}a_nz^n$, then by elementary computations, we
see that  $a_3=-1/6$ in both cases $\omega(z)=z$ and $\omega(z)=-z$.
Therefore, the corresponding harmonic mappings $f_{25}(z)$ and $f_{26}(z)$
do not belong to $\mathcal {S}_H^0(\frac{1}{2}\mathbb{Z})$.


\bsca The case $\varphi(z)=z(2-z)/(2(1-z)^2)$ and $\omega(z)=\pm z$. \esca
If $\omega(z)=z$, then
$$h'(z)=\frac{1}{(1-z)^4}, ~\mbox{ i.e. }~ h(z)=\frac{1}{3}\left ( \frac{1}{(1-z)^3}-1\right ),
$$
and if we let $h(z)=\sum_{n=1}^{\infty}a_nz^n$, then we see that $a_3=10/3$.

Similarly, for the case $\omega(z)=-z$, we obtain
$$h'(z)=\frac{1}{(1-z)^3(1+z)}=\frac{1}{2(1-z)^3}+\frac{1}{4(1-z)^2}+\frac{1}{4(1-z^2)}
$$
and the integration gives
\beqq
h(z)&=&\frac{1}{4(1-z)^2}+\frac{1}{4(1-z)}-\frac{1}{2}+\frac{1}{8}\log\left (\frac{1+z}{1-z}\right )\\
&=&\sum_{n=1}^{\infty}\left (\frac{n+2}{4}+\frac{1+(-1)^{n+1}}{8n}\right )z^n.
\eeqq
Clearly,  the corresponding harmonic mappings $f_{27}(z)$ and $f_{28}(z)$
do not belong to $\mathcal {S}_H^0(\frac{1}{2}\mathbb{Z})$.


\bsca The case $\varphi(z)=z(2+z)/(2(1+z)^2)$ and $\omega(z)=\pm z$. \esca
As in the previous subcase, for $\omega(z)=z$, we find that
$$h'(z)=\frac{1}{(1+z)^3(1-z)},
$$
and therefore,
\beqq
h(z) &=&-\frac{1}{4(1+z)^2}-\frac{1}{4(1+z)}+\frac{1}{2} +\frac{1}{8}\log\left (\frac{1+z}{1-z}\right )\\
&=&\sum_{n=1}^{\infty}(-1)^{n+1}\left (\frac{n+2}{4}+\frac{1+(-1)^{n+1}}{8n}\right )z^n.
\eeqq

Similarly, for $\omega(z)=-z$, it is easy to check that
$$h'(z)=\frac{1}{(1+z)^4}, ~\mbox{ i.e. }~ h(z)=-\frac{1}{3}\left ( \frac{1}{(1-z)^3}-1\right ),
$$
and for $h(z)=\sum_{n=1}^{\infty}a_nz^n$, we see that $a_3=10/3$.

Thus,  the corresponding harmonic mappings $f_{29}(z)$ and $f_{30}(z)$
do not belong to $\mathcal {S}_H^0(\frac{1}{2}\mathbb{Z})$.
\hfill $\Box$

\br
For a given $\theta \in [0,2\pi)$, let ${\mathcal M}(\theta)$ denote the set of all harmonic
functions $f=h+\overline{g}$ in $\ID$, with $h'(0)-1=0=h(0),$ that satisfies
$$g'(z)=e^{i\theta}zh'(z) ~\mbox{ and }~
{\rm Re}\left(1+z\frac{h''(z)}{h'(z)}\right)>-\frac{1}{2} ~\mbox{ for all $z\in \ID$}.
$$
In \cite[Theorem 1]{Bhara-samy-pre11}, the authors have shown that (without the normalization restriction)
every function in ${\mathcal M}(\theta)$ is univalent and close-to-convex in $\ID$. The case $\theta =0$ was
originally solved by Bshouty and Lyzzaik  \cite{Bshouty-Lyzzaik-2010}.



Note that if $f=h+\overline{g}\in {\mathcal M}(\pi)$ and $h(z)=\sum_{n=1}^{\infty}a_nz^n$ ($a_1=1$),
then it is a simple exercise to see that the co-analytic part $g$ of $f$ has the form
$$g(z)=-\sum_{n=2}^{\infty}a_{n-1}\left (\frac{n-1}{n}\right )z^n.
$$
In particular, the function
$$f_9(z)= {\rm Re}\left ( \frac{z}{1-z} \right ) +i {\rm Im}\left ( \frac{z}{(1-z)^2} \right )=
\sum_{n=1}^{\infty}\left (\frac{n+1}{2}\right ) z^n -
\overline{\sum_{n=2}^{\infty}\left (\frac{n-1}{2}\right ) z^n}
$$
belongs to ${\mathcal M}(\pi)$. From Theorem \ref{thm4.1}, it follows that
$f_9\in \mathcal {S}_H^0(\frac{1}{2}\IZ)$ and is convex in real direction. Furthermore,
the function $f_9(z)$ is a well-known extremal function for the convex
family of functions from  ${\mathcal S}_H^0$, and has the dilation $\omega (z)=-z$, see \cite{CS}.
This function has a special role in deriving convolutions results for harmonic mappings,
see \cite{Do,DN,LiPo1}. Thus, it is natural to ask if the convolution results of
\cite{Do} can be derived by using the class ${\mathcal M}(\theta)$. For example, when can
the harmonic convolution $f_1\ast f_2$ of functions in $f_1\in {\mathcal M}(\theta _1)$ and
$f_2\in {\mathcal M}(\theta _2)$ be univalent in $\ID$? We ask whether
${\mathcal M}(\theta)$ is included in the class of functions convex in a direction.

In \cite{Bhara-samy-pre11}, it was conjectured that the functions in ${\mathcal M}(0)$ are starlike in $\ID$.
The function $f_3(z)$ given by (see Subcase 2)
$$f_3(z)= {\rm Re}\left ( \frac{z}{(1-z)^2} \right ) +i {\rm Im}\left ( \frac{z}{1-z} \right )
$$
belongs to ${\mathcal M}(0)$. We recall that $f_3\in \mathcal {S}_H^0(\frac{1}{2}\IZ)$ and is convex
in real direction. However, the graph of the function $f_3(z)$ shows that $f_3(z)$ is not starlike
in $\ID$ (see the relevant figure in Section \ref{sec5}). This fact may be easily verified if one
uses for example, the method of proof of
Clunie and Sheil-Small \cite{CS} via the transformation
\[\zeta = \frac{1+z}{1-z}=\xi + i\eta \quad (\xi >0),
\]
and then writing
$$f_3(z)= \frac{1}{4}{\rm Re}\bigg[\bigg(\frac{1+z}{1-z}\bigg)^2-1\bigg]+i\frac{1}{2}{\rm Im}\bigg(\frac{1+z}{1-z}-1\bigg).
$$
An equivalent requirement for starlikness of an univalent harmonic mappings $f$ is
that $\arg f(e^{it})$ be nondecreasing function of $t$, i.e.,
$$ \frac{\partial }{\partial t} \arg  f(e^{it}) = {\rm Re} \left ( \frac{Df(e^{it})}{f(e^{it})}\right )  \geq 0,
$$
where $Df(z)=zf_{z}(z)-\overline{z}f_{\overline{z}}(z)$. In our case, a simple calculation shows that for $|t|<\pi/2$,
$$ \frac{\partial }{\partial t} \arg  f_3(e^{it}) = \frac{2\cos t}{-3+\cos 2t}  <0,
$$
showing that $f_3(z)$ is not starlike with respect to the origin.
Thus, there exists a function in ${\mathcal M}(0)$ that is not starlike in $\ID$. This example shows that
the conjecture is false.
\hfill $\Box$\er


\vspace{0.5cm}

\subsection{Proof of Theorem \ref{thm4.2}}
Let $f=h+\overline{g}\in \mathcal {S}_{H}^0(\frac{1}{2}\mathbb{Z})$ such that $f\in \mathcal
{C}\mathcal {V}(i)$, and further let $\psi=h+g$, where
$$h(z)=z+\sum_{n=2}^{\infty}a_n z^n~\mbox{ and }~g(z)=\sum_{n=2}^{\infty}b_n z^n.
$$
By Lemma \Ref{lem2.4} (with $\alpha =\pi/2$),  since $f$ is univalent, $f\in \mathcal {C}\mathcal {V}(i)$ if and
only if  $\psi $ is a univalent mapping such that $\psi\in \mathcal {C}\mathcal {V}(i)$. Moreover,
$\psi$ has half-integer coefficients, since $f\in \mathcal {S}_{H}^0(\frac{1}{2}\mathbb{Z})$.
It follows from Lemma \Ref{lem4.2} that $\psi\in \mathcal {T}_5$, where
$\mathcal {T}_{5}$ is given by \eqref{eq-half-shear1b}.

As in the proof of Theorem \ref{thm3.1}, since $f$ is sense-preserving, it follows that
$$-\frac{g'(z)}{\psi'(z)}\prec \frac{z}{1-z}.
$$
By Lemma \Ref{lem2.6}, we have $|2b_2|\leq 1$. Thus, we must have $b_2=0$ or $b_2=\pm 1/2$,
because $2b_2$ is an integer. Again, we observe that these two restrictions on $b_2$
are only necessary conditions for $f\in \mathcal {S}_{H}^0(\frac{1}{2}\mathbb{Z})$.

\bca The case $b_2=0$. \eca
As in the proof of Case 1 in the proof of Theorem \ref{thm4.1}, we can easily conclude that
 $g(z)\equiv0$. Hence $f$ is one of the functions from the set $\mathcal{T}_5$.

\medskip

Since $b_2=\pm 1/2$ are necessary, but not sufficient, for $f$ to be in $\mathcal {S}_H^0(\frac{1}{2}\mathbb{Z})$,
we need to deal these two cases separately.

\bca The case $b_2=\pm 1/2$.\eca

Since $b_2=\pm 1/2$, we can easily deduce from Lemma \Ref{lem2.5} that $\omega(z)=\pm z$.
Solving $g'(z)=\omega(z)h'(z)$ and $\psi(z)=h(z)+g(z)$, we obtain
\be\label{eq4.1.1}
h'(z)=\frac{\psi'(z)}{1+\omega(z)},
\ee
which gives us the harmonic function
\be\label{eq4.1.2}
f(z)=2i{\rm Im }h(z)+\overline{\psi(z)},
~\mbox{ with }~ h(z)=\int_{0}^z \frac{\psi'(t)}{1+\omega(t)}dt,
\ee
where $\omega(z)=\pm z$. We divide the remaining part of the proof into several subcases.

\setcounter{subcase}{0}

\bsca The case $\psi(z)=z$ and $\omega(z)=\pm z$.\esca

Evaluating the integral in \eqref{eq4.1.2} with $\psi(z)=\pm z$
yields
$$h(z)=
\left\{ \begin{array}{rl}
\ds\log(1+z)  & \mbox{if } \omega (z)=z,\\
\ds -\log(1-z) & \mbox{if } \omega(z)=-z.
\end{array}\right.
$$
Thus, the corresponding functions in $\mathcal{CV}(i)$ in these two cases are given by
$$ f_{1}(z)=\log(1+z)+\overline{z-\log(1+z)}=2i \arg (1+z) +\overline{z}
$$
and
$$f_{2}(z)=-\log(1-z)+\overline{z+\log(1-z)}=-2i \arg (1-z) +\overline{z},
$$
respectively. Clearly,  $f_1(z)$ and $f_2(z)$ do not belong to $\mathcal {S}_{H}^0(\frac{1}{2}\mathbb{Z})$.


\bsca The case  $\psi(z)=z/(1-z)$ and $\omega(z)=\pm z$. \esca

Evaluating the integral in \eqref{eq4.1.2} with
$\psi(z)=z/(1-z)$ and $\omega(z)=z$, we obtain
$$h(z)= \frac{z}{2(1-z)}+\frac{1}{4}\log\left (\frac{1+z}{1+z}\right )
=\sum_{n=1}^{\infty}\left(\frac{1}{2}+\frac{1+(-1)^{n+1}}{4n}\right )z^n.
$$
Clearly, the corresponding function $f_3(z)$ belongs to $\mathcal{CV}(i)$ but
do not have half-integer coefficients.

On the other hand, when $\psi(z)=z/(1-z)$ and $\omega(z)=-z$, it yields
$$ h(z)=\frac{1}{2(1-z)^2}-\frac{1}{2}=\frac{1}{2}\left ( \frac{z}{(1-z)^2}+
\frac{z}{1-z}\right )=
\sum_{n=1}^{\infty}\frac{n+1}{2}z^n
$$
so that
$$g(z)=\psi(z) -h(z) =-\frac{1}{2}\left ( \frac{z}{(1-z)^2}-\frac{z}{1-z}\right )=
-\sum_{n=2}^{\infty}\frac{n-1}{2}z^n.
$$
Consequently,  the corresponding harmonic function $f_4(z)$ in $\mathcal{CV}(i)$ takes the form
$$f_4(z)= {\rm Re}\left ( \frac{z}{1-z} \right ) +i {\rm Im}\left ( \frac{z}{(1-z)^2} \right ),
$$
which clearly belongs to $\mathcal {S}_{H}^0(\frac{1}{2}\mathbb{Z})$.

\bsca The case $\psi(z)=z/(1+z)$ and $\omega(z)=\pm z$. \esca

In this case, it follows from \eqref{eq4.1.2} that
$$h(z)=\left\{ \begin{array}{rl}
\ds -\frac{1}{2(1+z)^2}+\frac{1}{2}  & \mbox{if } \omega (z)=z,\\
\ds \frac{z}{2(1+z)}+\frac{1}{4}\log\left( \frac{1+z}{1-z}\right )  & \mbox{if } \omega (z)=-z.
\end{array}\right.
$$
In the case $\omega (z)=z$, $h(z)$ above may be rewritten as
$$h(z)=\frac{1}{2}\left ( \frac{z}{(1+z)^2}+
\frac{z}{1+z}\right )=\sum_{n=1}^{\infty}(-1)^{n+1}\frac{n+1}{2}z^n
$$
and therefore, the corresponding co-analytic part is
$$g(z)=\psi(z) -h(z) =-\frac{1}{2}\left ( \frac{z}{(1+z)^2}-\frac{z}{1+z}\right )=
-\sum_{n=2}^{\infty}(-1)^{n+1}\frac{n-1}{2}z^n.
$$
Thus, the corresponding harmonic function $f_5(z)$ in $\mathcal{CV}(i)$ takes the form
$$f_5(z)= {\rm Re}\left ( \frac{z}{1+z} \right ) +i {\rm Im}\left ( \frac{z}{(1+z)^2} \right )
$$
and is clearly in $\mathcal {S}_{H}^0(\frac{1}{2}\mathbb{Z})$.

In the case $\omega (z)=-z$, it is again easy to see that the corresponding harmonic function
$f_6(z)$ in $\mathcal{CV}(i)$ does not have half-integer coefficients, since the analytic part of
$f_6(z)$ given by
$$h(z)= \frac{z}{2(1+z)}+\frac{1}{4}\log\left( \frac{1+z}{1-z}\right )
=\sum_{n=1}^{\infty}(-1)^{n+1}\left (\frac{1}{2}+\frac{1+(-1)^{n+1}}{4n}\right )z^n
$$
does not have half-integer coefficients.

\bsca The case $\psi(z)=z/(1-z^2)$ and $\omega(z)=\pm z$.
\esca

Evaluating $h'(z)$ with $\psi(z)=\frac{z}{1-z^2}$ and $\omega(z)=\pm
z$ in \eqref{eq4.1.1} yields

$$
h'(z)=\left\{ \begin{array}{rl}
\ds\frac{1+z^2}{(1-z^2)^2(1+z)}& \mbox{for } \omega(z)=z,\\
\ds \frac{1+z^2}{(1-z^2)^2(1-z)}& \mbox{for } \omega(z)=-z.
\end{array}\right.
$$
It is easy to check that in these two cases, the third coefficient of $h(z)$ is $a_3=4/3$, and
so, the corresponding harmonic functions $f_7(z)$ and  $f_8(z)$
do not have half-integer coefficients.


\bsca The case $\psi(z)= z(2-z)/(2(1-z))$ and $\omega(z)=\pm
z$. \esca

It follows from the integral in \eqref{eq4.1.2} with
$\psi(z)=z(2-z)/(2(1-z))$ and $\omega(z)= z$ that
$$h(z)=\frac{5}{8}\log(1+z)-\frac{1}{8}\log(1-z)+\frac{z}{4(1-z)}
=\sum_{n=1}^{\infty}\left (\frac{5(-1)^{n+1}+1}{8n}+\frac{1}{4}\right )z^n,
$$
which implies that the corresponding harmonic function $f_9(z)$ does not belong to
$\mathcal {S}_{H}^0(\frac{1}{2}\mathbb{Z})$.

In the case $\omega(z)=- z$, it yields
$$h(z)=\frac{1}{4(1-z)^2}-\frac{1}{2}\log(1-z)-\frac{1}{4}
=\sum_{n=1}^{\infty}\left (\frac{n+1}{4}+\frac{1}{2n}\right )z^n,
$$
and hence, the corresponding harmonic function $f_{10}(z)$ does not have half-integer coefficients.


\bsca The case $\psi(z)= z(2+z)/(2(1+z))$ and $\omega(z)=\pm z$. \esca
From the integral in \eqref{eq4.1.2} with
$\psi(z)=\frac{z(2+z)}{2(1+z)}$ and $\omega(z)= z$, we obtain
\beqq
h(z)&=&\frac{1}{8}\log(1+z)-\frac{5}{8}\log(1-z) +\frac{z}{4(1+z)}\\
&=&\sum_{n=1}^{\infty}(-1)^{n+1}\left (\frac{5(-1)^{n+1} +1}{8n}+\frac{1}{4}\right )z^n.
\eeqq

In the case $\omega(z)=-z$, we get
$$h(z)=-\frac{1}{4(1+z)^2}+\frac{1}{2}\log(1+z)+\frac{1}{4}
=\sum_{n=1}^{\infty}(-1)^{n+1}\left (\frac{n+1}{4}+\frac{1}{2n}\right )z^n.
$$

Hence in these two cases, the corresponding harmonic functions $f_{11}(z)$
and $f_{12}(z)$ do not belong to ${S}_{H}^0(\frac{1}{2}\mathbb{Z})$.


\bsca The case $\psi(z)=z(2-z^2)/(2(1-z^2))$ and $\omega(z)=\pm z$.
\esca

When $\omega(z)=z$, it follows from \eqref{eq4.1.1} that
$$h'(z)=\frac{2-z^2+z^4}{2(1-z^2)^2(1+z)}
=\frac{1}{4(1+z)^3}+\frac{1}{8(1-z)^2}+\frac{1}{16(1-z)}+\frac{9}{16(1+z)},$$
and therefore,
\beqq
h(z)&=&-\frac{1}{8(1+z)^2}+\frac{1}{8(1-z)}-\frac{1}{16}\log(1-z)+\frac{9}{16}\log(1+z)\\
&=&\sum_{n=1}^{\infty}\left (\frac{(n+1)(-1)^{n+1} +1}{8} +\frac{1+9(-1)^{n+1}}{16n}\right )z^n.
\eeqq

If  $\omega(z)=-z$, then from \eqref{eq4.1.1} we get
$$h'(z)=\frac{2-z^2+z^4}{2(1-z^2)^2(1-z)}=\frac{1}{4(1-z)^3}+\frac{1}{8(1+z)^2}+\frac{1}{16(1+z)}+\frac{9}{16(1-z)},
$$
and thus
\beqq
h(z)&=&\frac{1}{8(1-z)^2}-\frac{1}{8(1+z)}+\frac{1}{16}\log(1+z)-\frac{9}{16}\log(1-z)\\
&=&\sum_{n=1}^{\infty}(-1)^{n+1}\left (\frac{(n+1)(-1)^{n+1}+1}{8}
+\frac{1+9(-1)^{n+1}}{16n}\right )z^n.
\eeqq
In these two cases, the corresponding harmonic functions $f_{13}(z)$
and $f_{14}(z)$ do not have half-integer coefficients.


\bsca The case $\psi(z)=z(2- z)/(2(1-z^2))$ and $\omega(z)=\pm z$. \esca

Evaluating $h'(z)$ and $h(z)$ in \eqref{eq4.1.1} and \eqref{eq4.1.2} with
$\psi(z)=\frac{z(2- z)}{2(1-z^2)}$ and $\omega(z)=z$, we obtain
$$ h'(z)=\frac{1-z+z^2}{(1-z^2)^2(1+z)}=\frac{3}{4(1+z)^3}
+\frac{1}{8(1-z)^2}+\frac{1}{16}\left (\frac{1}{1-z}+\frac{1}{1+z}\right )
$$
and the integration gives
\beqq
h(z)&=&\frac{1}{16}\log\left (\frac{1+z}{1-z}\right )+\frac{1}{8(1-z)}-\frac{3}{8(1+z)^2}+\frac{1}{4}\\
&=& \sum_{n=1}^{\infty}\left (\frac{(-1)^{n+1}+1}{16n}
+\frac{1+3(n+1)(-1)^{n+1}}{8}\right )z^n.
\eeqq

If $\psi(z)=\frac{z(2- z)}{2(1-z^2)}$ and $\omega(z)=-z$, it follows that
$$h'(z)=\frac{1-z+z^2}{(1-z^2)^2(1-z)}=\frac{1}{4(1-z)^3}+\frac{3}{8(1+z)^2}+\frac{3}{16(1-z)}+\frac{3}{16(1+z)},
$$
and thus integration gives
\beqq
h(z)&=&\frac{3}{16}\log\left (\frac{1+z}{1-z}\right )+\frac{1}{8(1-z)^2}-\frac{3}{8(1+z)}+\frac{1}{4}\\
&=& \sum_{n=1}^{\infty}\left (\frac{3(-1)^{n+1}+3}{16n}+\frac{n+1+3(-1)^{n+1}}{8}\right )z^n.
\eeqq

Therefore, in these two cases, the corresponding harmonic functions $f_{15}(z)$
and $f_{16}(z)$ do not belong to $\mathcal {S}_{H}^0(\frac{1}{2}\mathbb{Z})$.


\bsca The case $\psi(z)= z(2+z)/(2(1-z^2))$ and $\omega(z)=\pm z$.
\esca
If $\psi(z)=\frac{z(2+z)}{2(1-z^2)}$ and $\omega(z)=z$,
then it follows from \eqref{eq4.1.1} and \eqref{eq4.1.2} that
$$ h'(z)=\frac{1+z+z^2}{(1-z^2)^2(1+z)}=\frac{1}{4(1+z)^3}+\frac{3}{8(1-z)^2}+\frac{3}{16(1+z)}+\frac{3}{16(1-z)},
$$
and therefore,
\beqq
h(z)&=&\frac{3}{16}\log\left (\frac{1+z}{1-z}\right )-\frac{1}{8(1+z)^2}+\frac{3}{8(1-z)}-\frac{1}{4}\\
&=& \sum_{n=1}^{\infty}(-1)^{n+1}\left (\frac{3(-1)^{n+1} +3}{16n}+\frac{ n+1+ 3(-1)^{n+1}}{8}\right )z^n.
\eeqq

When $\psi(z)=\frac{z(2+z)}{2(1-z^2)}$ and $\omega(z)=-z$, we have
$$
h'(z)=\frac{1+z+z^2}{(1-z^2)^2(1-z)}=\frac{3}{4(1-z)^3}+\frac{1}{8(1+z)^2}+\frac{1}{16(1+z)}+\frac{1}{16(1-z)},
$$
and therefore,
\beqq
h(z)&=&\frac{1}{16}\log\left (\frac{1+z}{1-z}\right )
+\frac{3}{8(1-z)^2}-\frac{1}{8(1+z)}-\frac{1}{4}\\
&=&\sum_{n=1}^{\infty}(-1)^{n+1}\left (\frac{(-1)^{n+1} +1}{16n}
+\frac{1+3(n+1)(-1)^{n+1}}{8}\right )z^n.
\eeqq
Thus, in these two cases, the corresponding harmonic functions $f_{17}(z)$
and $f_{18}(z)$ do not belong to $\mathcal {S}_{H}^0(\frac{1}{2}\mathbb{Z})$.
\hfill $\Box$

\br
In the proof of Theorem \ref{thm4.1}, we observe that there are thirty functions that are convex
in real direction (for the case $b_2=\pm 1/2$) out of which only six have half-integer coefficients.

From the proof of Theorem \ref{thm4.1}, we see that there are eighteen functions that are convex
in vertical direction (for the case $b_2=\pm 1/2$) but only two of these have half-integer coefficients.
\er

\subsection*{Acknowledgements}
The work of Ms. Jinjing Qiao was supported by Centre for International Co-operation in Science (CICS)
through the award of ``INSA JRD-TATA Fellowship." The work was completed
during her visit in April-June 2012 to the Department of Mathematics, IIT Madras, Chennai.
The author thanks both CICS and IIT Madras for their support.



\begin{thebibliography}{17}

\bibitem{Bhara-samy-pre11}
S.~V.~Bharanedhar and S.~Ponnusamy,
\textrm{Coefficient conditions for harmonic univalent mappings and hypergeometric mappings},
\textit{Rocky Mountain J. Math.} (2012), To appear.

\bibitem{Bi}L. Bieberbach,
\"{U}ber die Koeffizienten derjenigen Potenzreihen, welche eine schlighte
Abbildung des Einheitskreis vermitteln,
\textit{S.-B. Preuss. Akad. Wiss.} \textbf{38}(1916), 940--955.

\bibitem{Bshouty-Lyzzaik-2010} D.~Bshouty and A.~Lyzzaik,
\textrm{Close-to-convexity criteria for planar harmonic mappings},
\textit{Complex Anal. Oper. Theory} {\bf 5}(3)(2011), 767--774.

\bibitem{CS} J. G. Clunie and T. Sheil-Small,
Harmonic univalent functions,
\textit{Ann. Acad. Sci. Fenn. Ser. A. I.} \textbf{9}(1984), 3--25.

\bibitem{Do} M.~Dorff,
\textrm{Convolutions of planar harmonic convex mappings},
\textit{Comp. Vari. Theo. Appl.} {\bf 45}(2001), 263--271.

\bibitem{DN} M.~Dorff, M.~Nowak and M.~Wo{\l}oszkiewicz,
\textrm{Convolutions of harmonic convex mappings,}
\textit{Comp. Vari. Elliptic Eqn.} {\bf 57}(5)(2012), 489--503.

\bibitem{Du} P. Duren,
Univalent Functions,
\textit{Springer-Verlag}, New York, Berlin, Heidelberg, Tokyo, 1982.

\bibitem{Du1} P. Duren,
Harmonic mappings in the plane,
\textit{Cambridge Univ. Press}, 2004.

\bibitem{Fr} B. Friedman,
Two theorems on schlicht functions,
\textit{Duke Math. J.} \textbf{13}(1946), 171--177.

\bibitem{Go} A. W. Goodman,
Univalent functions, Vols. 1-2,
Mariner, Tampa, Florida, 1983.

\bibitem{Gre} P. Greiner,
Geometric properties of harmonic shears,
\textit{Comput. Methods Funct. Theory} {\bf 4}(1)(2004), 77--96.

\bibitem{Gr} T. H. Gronwall,
Some remarks on conformal representation,
\textit{Ann. of Math.} \textbf{16}(1914-1915), 72--76.

\bibitem{HengSch70} W. Hengartner and G. Schober,
On schlicht mappings to domains convex in one direction,
\textit{Comment. Math. Helv.} \textbf{45}(1970), 303-–314.

\bibitem{HS} N. Hiranuma and T. Sugawa,
Univalent functions with half-integer coefficients, Preprint.

\bibitem{Jen87} J. A. Jenkins, On univalent functions with integral coefficients,
\textit{Complex Variables Theory Appl.} \textbf{9}(1987), 221--226.

\bibitem{Lecko02} A. Lecko,
On the class of functions convex in the negative direction of the imaginary axis,
\textit{J. Aust. Math. Soc.} \textbf{73}(2002), 1--10.

\bibitem{Le} H. Lewy,
On the non-vanishing of the Jacobian in certain one-to-one mappings,
\textit{Bull. Amer. Math. Soc.} \textbf{42}(1936), 689--692.

\bibitem{Li55} V. Linus,
Note on univalent functions,
\textit{Amer. Math. Monthly} \textbf{62}(1955), 109--110.

\bibitem{LiPo1} Liulan Li and S.~Ponnusamy, Solution to an open problem on convolutions
of harmonic mappings,
\textit{Comp. Vari. Elliptic Eqn.} (2012), Accepted.


\bibitem{OP} M. Obradovi\'{c} and S. Ponnusamy,
New criteria and distortion theorems for univalent functions,
\textit{Complex Variables Theory Appl.} \textbf{44}(2001), 173--191.


\bibitem{P} Ch. Pommerenke,
Univalent functions, Vandenhoeck and
Ruprecht, G\"ottingen, 1975.

\bibitem{Rober36} M. S. Robertson,
Analytic functions starlike in one direction,
\textit{Am. J. Math.} \textbf{58}(1936), 465--472.

\bibitem{Rog} W. W. Rogosinski,
On the coefficients of subordinate functions,
\textit{Proc. London Math. Soc.} \textbf{48}(1943), 48--82.

\bibitem{Roy55} W. C. Royster,
Rational univalent functions,
\textit{Amer. Math. Monthly} \textbf{63}(1956), 326--328.

\bibitem{RZ} W. C. Royster and M. Ziegler,
Univalent functions convex in one direction,
\textit{Publ. Math. Debrecen} \textbf{23}(1976), 339--345.

\bibitem{Sch2001} L. E. Schaubroeck,
Growth, distortion and coefficient bounds for plane harmonic mappings
convex in one direction,
\textit{Rocky Mountain J. Math.} \textbf{31}(2)(2001), 624--639

\bibitem{Sh} T. -S. Shah,
 On the coefficients of schlicht functions,
\textit{J. Chinese Math. Soc. (N. S.)} \textbf{1}(1951), 98--107.

\bibitem{To} S. B. Townes,
A theorem on schlicht functions,
\textit{Proc. Amer. Math. Soc.} \textbf{5}(1954), 585--588.

\end{thebibliography}
\end{document}